\documentclass{amsart}
\usepackage{amssymb}
\newtheorem{theorem}{Theorem}

\newtheorem{lemma}[theorem]{Lemma}

\newtheorem{question}[theorem]{Question}

\newtheorem{definition}[theorem]{Definition}

\newcommand{\QED}{\end{proof}}
\newenvironment{proclamation}{\smallskip\noindent\proclaim}{\par\smallskip}

\def\proclaim[#1]{{\bf #1}}
\def\BF#1.{{\bf #1.}}
\newcommand{\url}[1]{{\tt #1}}

\newcommand{\Godel}{G\"odel}

\newcommand{\B}{{\mathbb B}}
\newcommand{\C}{{\mathbb C}}

\newcommand{\N}{{\mathbb N}}
\renewcommand{\P}{{\mathbb P}}
\newcommand{\Q}{{\mathbb Q}}

\newcommand{\R}{{\mathbb R}}

\newcommand{\Vbar}{{\overline{V}}}

\newcommand{\of}{\subseteq}

\newcommand{\set}[1]{\{\,{#1}\,\}}

\newcommand{\Con}{\mathop{{\rm Con}}}

\newcommand{\satisfies}{\models}
\newcommand{\forces}{\Vdash}
\newcommand{\proves}{\vdash}
\newcommand{\possible}{\mathop{\raisebox{-1pt}{$\Diamond$}}}
\newcommand{\necessary}{\mathop{\raisebox{-1pt}{$\Box$}}}

\newcommand{\axiomf}[1]{{\rm #1}}
\newcommand{\theoryf}[1]{\hbox{$\mathsf{#1}$}}

\newcommand{\Mantle}{{\mathord{\rm M}}}
\newcommand{\gMantle}{{\mathord{\rm gM}}}  % the generic mantle
\newcommand{\gHOD}{\ensuremath{\mathord{{\rm g}\HOD}}} % the limit HOD

\newcommand{\Union}{\bigcup}

\newcommand{\intersect}{\cap}
\newcommand{\Intersect}{\bigcap}

\newcommand{\smalllt}{\mathrel{\mathchoice{\raise2pt\hbox{$\scriptstyle<$}}{\raise1pt\hbox{$\scriptstyle<$}}{\raise0pt\hbox{$\scriptscriptstyle<$}}{\scriptscriptstyle<}}}
\newcommand{\smallleq}{\mathrel{\mathchoice{\raise2pt\hbox{$\scriptstyle\leq$}}{\raise1pt\hbox{$\scriptstyle\leq$}}{\raise1pt\hbox{$\scriptscriptstyle\leq$}}{\scriptscriptstyle\leq}}}
\newcommand{\lt}{\smalllt}

\newcommand{\boolval}[1]{\mathopen{\lbrack\!\lbrack}\,#1\,\mathclose{\rbrack\!\rbrack}}
\def\[#1]{\boolval{#1}}

\newcommand{\UnderTilde}[1]{{\setbox1=\hbox{$#1$}\baselineskip=0pt\vtop{\hbox{$#1$}\hbox to\wd1{\hfil$\sim$\hfil}}}{}}
\newcommand{\Undertilde}[1]{{\setbox1=\hbox{$#1$}\baselineskip=0pt\vtop{\hbox{$#1$}\hbox to\wd1{\hfil$\scriptstyle\sim$\hfil}}}{}}
\newcommand{\undertilde}[1]{{\setbox1=\hbox{$#1$}\baselineskip=0pt\vtop{\hbox{$#1$}\hbox to\wd1{\hfil$\scriptscriptstyle\sim$\hfil}}}{}}
\newcommand{\UnderdTilde}[1]{{\setbox1=\hbox{$#1$}\baselineskip=0pt\vtop{\hbox{$#1$}\hbox to\wd1{\hfil$\approx$\hfil}}}{}}
\newcommand{\Underdtilde}[1]{{\setbox1=\hbox{$#1$}\baselineskip=0pt\vtop{\hbox{$#1$}\hbox to\wd1{\hfil\scriptsize$\approx$\hfil}}}{}}

\newcommand{\st}{\mid}

\def\<#1>{\langle#1\rangle}

\newcommand{\val}{\mathop{\rm val}\nolimits}

\newcommand{\ORD}{\mathop{{\rm ORD}}}
\newcommand{\REG}{\mathop{{\rm REG}}}

\newcommand{\ZFC}{{\rm ZFC}}
\newcommand{\ZF}{{\rm ZF}}

\newcommand{\CH}{{\rm CH}}
\newcommand{\KP}{{\rm KP}}

\newcommand{\AC}{{\rm AC}}

\newcommand{\AFA}{{\rm AFA}}
\newcommand{\AS}{{\rm AS}}

\newcommand{\GA}{{\rm GA}}
\newcommand{\DDG}{{\rm DDG}}

\newcommand{\MM}{{\rm MM}}

\newcommand{\PFA}{{\rm PFA}}

\newcommand{\HOD}{{\rm HOD}}

\newcommand{\PA}{{\rm PA}}

\newcommand{\cell}[1]{\boxit{\hbox to 17pt{\strut\hfil$#1$\hfil}}}
\newcommand{\head}[2]{\lower2pt\vbox{\hbox{\strut\footnotesize\it\hskip3pt#2}\boxit{\cell#1}}}
\newcommand{\boxit}[1]{\setbox4=\hbox{\kern2pt#1\kern2pt}\hbox{\vrule\vbox{\hrule\kern2pt\box4\kern2pt\hrule}\vrule}}
\newcommand{\Col}[3]{\hbox{\vbox{\baselineskip=0pt\parskip=0pt\cell#1\cell#2\cell#3}}}
\newcommand{\tapenames}{\raise 5pt\vbox to .7in{\hbox to .8in{\it\hfill input: \strut}\vfill\hbox to
.8in{\it\hfill scratch: \strut}\vfill\hbox to .8in{\it\hfill output: \strut}}}
\newcommand{\Head}[4]{\lower2pt\vbox{\hbox to25pt{\strut\footnotesize\it\hfill#4\hfill}\boxit{\Col#1#2#3}}}
\newcommand{\Dots}{\raise 5pt\vbox to .7in{\hbox{\ $\cdots$\strut}\vfill\hbox{\ $\cdots$\strut}\vfill\hbox{\
$\cdots$\strut}}}
\newcommand{\df}{\it} % use italic for definition terms. Idea: also use this to create an index of definitions, if MakeIndex is true.
\usepackage{diagrams}
\diagramstyle[tight,centredisplay%
,dpi=600%              for more modern laser-printers
,noPostScript%
,heads=LaTeX%
]%

\begin{document}
\author{Joel David Hamkins}
\address{J. D. Hamkins, Department of Philosophy, New York University, 5
Washington Place, New York, NY 10003 \& Mathematics, The
Graduate Center of The City University of New York, 365
Fifth Avenue, New York, NY 10016 \& Mathematics, The
College of Staten Island of CUNY, Staten Island, NY 10314}
\email{jhamkins@gc.cuny.edu, http://jdh.hamkins.org} \today
\thanks{This article is based on a talk I gave at the
conference Philosophy of Mathematics, held at New York
University in April, 2009. My research has been supported
in part by grants from the CUNY Research Foundation, the
National Science Foundation (DMS-0800762) and the Simons
Foundation, for which I am very grateful. I would like to
thank the anonymous referee, as well as Haim Gaifman,
Victoria Gitman, Peter Koellner, Graham Leach-Krouse, Toby
Meadows, Barbara Montero, Bob Solovay, Sean Walsh and W.
Hugh Woodin for extremely helpful comments and
discussions.}
%\subjclass{03Axx,03Exx}

\begin{abstract}
The multiverse view in set theory, introduced and argued
for in this article, is the view that there are many
distinct concepts of set, each instantiated in a
corresponding set-theoretic universe. The universe view, in
contrast, asserts that there is an absolute background set
concept, with a corresponding absolute set-theoretic
universe in which every set-theoretic question has a
definite answer. The multiverse position, I argue, explains
our experience with the enormous range of set-theoretic
possibilities, a phenomenon that challenges the universe
view. In particular, I argue that the continuum hypothesis
is settled on the multiverse view by our extensive
knowledge about how it behaves in the multiverse, and as a
result it can no longer be settled in the manner formerly
hoped for.
\end{abstract}

\title{The set-theoretic multiverse}\maketitle

\section{Introduction}

Set theorists commonly take their subject as constituting
an ontological foundation for the rest of mathematics, in
the sense that abstract mathematical objects can be
construed fundamentally as sets, and because of this, they
regard set theory as the domain of all mathematics.
Mathematical objects {\it are} sets for them, and being
precise in mathematics amounts to specifying an object in
set theory. These sets accumulate transfinitely to form the
cumulative universe of all sets, and the task of set theory
is to discover its fundamental truths.

The {\df universe view} is the commonly held philosophical
position that there is a unique absolute background concept
of set, instantiated in the corresponding absolute
set-theoretic universe, the cumulative universe of all
sets, in which every set-theoretic assertion has a definite
truth value. On this view, interesting set-theoretic
questions, such as the continuum hypothesis and others,
have definitive final answers. Adherents of the universe
view often point to the increasingly stable consequences of
the large cardinal hierarchy, particularly in the realm of
projective sets of reals with its attractive determinacy
and regularity features, as well as the
forcing-absoluteness properties for $L(\R)$, as evidence
that we are on the right track towards the final answers to
these set theoretical questions. Adherents of the universe
view therefore also commonly affirm these large cardinal
and regularity features in the absolute set-theoretic
universe. The pervasive independence phenomenon in set
theory is described on this view as a distraction, a side
discussion about provability rather than truth---about the
weakness of our theories in finding the truth, rather than
about the truth itself---for the independence of a
set-theoretic assertion from \ZFC\ tells us little about
whether it holds or not in the universe.

In this article, I shall argue for a contrary position, the
{\it multiverse view}, which holds that there are diverse
distinct concepts of set, each instantiated in a
corresponding set-theoretic universe, which exhibit diverse
set-theoretic truths. Each such universe exists
independently in the same Platonic sense that proponents of
the universe view regard their universe to exist. Many of
these universes have been already named and intensely
studied in set theory, such as the classical models $L$ and
$\HOD$, the increasingly sophisticated definable inner
models of large cardinals and the vast diversity of forcing
extensions, although it is better to understand these
descriptions as relative construction methods, since the
resulting universe described depends on the initial
universe in which the constructions are undertaken. Often
the clearest way to refer to a set concept is to describe
the universe of sets in which it is instantiated, and in
this article I shall simply identify a set concept with the
model of set theory to which it gives rise. By adopting a
particular concept of set, we in effect adopt that universe
as our current mathematical universe; we jump inside and
explore the nature of set theory offered by that universe.
In this sense, the multiverse view does not undermine the
claim that set theory serves an ontological foundation for
mathematics, since one expects to find all the familiar
classical mathematical objects and structures inside any
one of the universes in the multiverse, but rather it is
directed at the claim that there is a unique absolute
background concept of set, whose set-theoretic truths are
immutable.

In particular, I shall argue in section
\ref{Section:CaseStudyCH} that the question of the
continuum hypothesis is settled on the multiverse view by
our extensive, detailed knowledge of how it behaves in the
multiverse. As a result, I argue, the continuum hypothesis
can no longer be settled in the manner formerly hoped for,
namely, by the introduction of a new natural axiom
candidate that decides it. Such a dream solution template,
I argue, is impossible because of our extensive experience
in the \CH\ and $\neg\CH$ worlds.

%After discussing in section
%\ref{Section:ChallengeOfDiversePossibility} the challenge
%to the universe view posed by the enormous range of
%set-theoretic possibility, I discuss in section
%\ref{Section:OntologyOfForcing} the ontological status of
%forcing extensions of the universe. In section
%\ref{Section:AnalogyWithGeometry}, I describe the
%fundamental analogy between multiverse set theory and the
%situation of geometry. I address several categoricity
%arguments in section \ref{Section:Categoricity}, and
%explain in section \ref{Section:Adjudication} the sense in
%which the multiverse provides the right context for
%adjudicating internal differences within the universe view.
%In section \ref{Section:CaseStudyCH}, I argue that the
%continuum hypothesis is settled by our extensive knowledge
%of how it behaves in the multiverse, and as a result it can
%no longer be settled in the manner formerly hoped. I
%describe the multiverse attitude to $V=L$ in section
%\ref{Section:CaseStudyV=L}, and then embark on a fuller
%description of the multiverse vision with the multiverse
%axioms in section \ref{Section.multiverseAxioms}.

The multiverse view is one of higher-order
realism---Platonism about universes---and I defend it as a
realist position asserting actual existence of the
alternative set-theoretic universes into which our
mathematical tools have allowed us to glimpse. The
multiverse view, therefore, does not reduce via proof to a
brand of formalism. In particular, we may prefer some of
the universes in the multiverse to others, and there is no
obligation to consider them all as somehow equal.

The assertion that there are diverse concepts of set is a
meta-mathematical as opposed to a mathematical claim, and
one does not expect the properties of the multiverse to be
available when undertaking an internal construction within
a universe. That is, we do not expect to see the whole
multiverse from within any particular universe.
Nevertheless, set theory does have a remarkable ability to
refer internally to many alternative set concepts, as when
we consider definable inner models or various outer models
to which we have access. In this way, set theory allows us
to mount largely mathematical explorations of questions
having a deeply philosophical nature, perhaps providing
mathematical footholds for further philosophical inquiry.
In the appendix of this article I describe two examples of
such work, the modal logic of forcing and set-theoretic
geology, which investigate the features of the
set-theoretic universe in the context of all its forcing
extensions and grounds.

\section{The challenge of diverse set-theoretic possibilities}
\label{Section:ChallengeOfDiversePossibility}

Imagine briefly that set theory had followed an alternative
history, that as the theory developed, theorems were
increasingly settled in the base theory; that the
independence phenomenon was limited to paradoxical-seeming
meta-logic statements; that the few true independence
results occurring were settled by missing natural
self-evident set principles; and that the basic structure
of the set-theoretic universe became increasingly stable
and agreed-upon. Such developments would have constituted
evidence for the universe view. But the actual history is
not like this.

Instead, the most prominent phenomenon in set theory has
been the discovery of a shocking diversity of set-theoretic
possibilities. Our most powerful set-theoretic tools, such
as forcing, ultrapowers and canonical inner models, are
most naturally and directly understood as methods of
constructing alternative set-theoretic universes.  A large
part of set theory over the past half-century has been
about constructing as many different models of set theory
as possible, often to exhibit precise features or to have
specific relationships with other models. Would you like to
live in a universe where $\CH$ holds, but $\Diamond$ fails?
Or where $2^{\aleph_n}=\aleph_{n+2}$ for every natural
number $n$? Would you like to have rigid Suslin trees?
Would you like every Aronszajn tree to be special? Do you
want a weakly compact cardinal $\kappa$ for which
$\Diamond_\kappa(\REG)$ fails? Set theorists build models
to order.

As a result, the fundamental objects of study in set theory
have become the models of set theory, and set theorists
move with agility from one model to another. While group
theorists study groups, ring theorists study rings and
topologists study topological spaces, set theorists study
the models of set theory. There is the constructible
universe $L$ and its forcing extensions $L[G]$ and
non-forcing extensions $L[0^\sharp]$; there are
increasingly sophisticated definable inner models with
large cardinals $L[\mu]$, $L[\vec E]$ and so on; there are
models $V$ with much larger large cardinals and
corresponding forcing extensions $V[G]$, ultrapowers $M$,
cut-off universes $L_\delta$, $V_\alpha$, $H_\kappa$,
universes $L(\R)$, $\HOD$, generic ultrapowers, boolean
ultrapowers and on and on and on. As for forcing
extensions, there are those obtained by adding one Cohen
real, or many, or by other c.c.c. or proper (or
semi-proper) forcing, or by long iterations of these, or by
the L\'evy collapse, or by the Laver preparation or by
self-encoding forcing, and on and on and on. Set theory
appears to have discovered an entire cosmos of
set-theoretic universes, revealing a category-theoretic
nature for the subject, in which the universes are
connected by the forcing relation or by large cardinal
embeddings in complex commutative diagrams, like
constellations filling a dark night sky.\footnote{This
situation was anticipated by Mostowski, who in a talk
recorded in
\cite{Lakatos1967:ProblemsInThePhilosophyOfMathematics}
asserted that there are ``several essentially different
notions of set which are equally admissible as an intuitive
basis for set theory,'' and Kalm\'ar agrees with him about
future practice (in 1965) by saying ``I guess that in the
future we shall say as naturally ``let us take a set
theory,'' as we take now a group G or a field F.''}

This abundance of set-theoretic possibilities poses a
serious difficulty for the universe view, for if one holds
that there is a single absolute background concept of set,
then one must explain or explain away as imaginary all of
the alternative universes that set-theorists seem to have
constructed. This seems a difficult task, for we have a
robust experience in those worlds, and they appear fully
set-theoretic to us. The multiverse view, in contrast,
explains this experience by embracing them as real, filling
out the vision hinted at in our mathematical experience,
that there is an abundance of set-theoretic worlds into
which our mathematical tools have allowed us to glimpse.

The case of the definable inner models, such as $L(\R)$ and
$\HOD$, may seem at first to be unproblematic for the
universe view, as they are directly accessible from their
containing universes. I counter this attitude, however, by
pointing out that much of our knowledge of these inner
models has actually arisen by considering them inside
various outer models. We understand the coquettish nature
of \HOD, for example, by observing it to embrace an entire
forcing extension, where sets have been made definable,
before relaxing again in a subsequent extension, where they
are no longer definable. In the case of the axiom of
choice, Cohen first moves to an outer model before
constructing the desired inner model of $\neg\AC$. In this
sense, the universe view lacks a full account of the
definable inner models.

But it is the outer models, of course, such as the diverse
forcing extensions, that most directly challenge the
universe view. We do have a measure of access into the
forcing extensions via names and the forcing relation,
allowing us to understand the objects and truths of the
forcing extension while remaining in the ground model. The
multiverse view explains our mathematical experience with
these models by positing that, indeed, these alternative
universes exist, just as they seem to exist, with a full
mathematical existence, fully as real as the universe under
the universe view. Thus, in our mathematical experience the
classical set concept splinters into a diverse array of
parallel set concepts, and there seems little reason to
think that we have discovered more than a tiny part of the
multiverse.

The multiverse view does not abandon the goal of using set
theory as an epistemological and ontological foundation for
mathematics, for we expect to find all our familiar
mathematical objects, such as the integer ring, the real
field and our favorite topological spaces, inside any one
of the universes of the multiverse. On the multiverse view,
set theory remains a foundation for the classical
mathematical enterprise. The difference is that when a
mathematical issue is revealed to have a set-theoretic
dependence, as in the independence results, then the
multiverse response is a careful explanation that the
mathematical fact of the matter depends on which concept of
set is used, and this is almost always a very interesting
situation, in which one may weigh the desirability of
various set-theoretic hypotheses with their mathematical
consequences. The universe view, in contrast, insists on
one true answer to the independent question, although we
seldom know which it is.

The beginning of a multiverse perspective appears already
with von Neumann in his paper ``An axiomatization of set
theory,'' (1925, available in \cite[p.
393-413]{vanHeijenoort1967:FromFregeToGodelSourcebookInMathLogic}),
where he considers (p. 412-413) the situation where one
model of set theory can be a set inside another model of
set theory and notes that a set that is ``finite'' in the
former may be infinite in the latter; and similarly, a
well-ordering in the former model may be ill-founded in the
latter. He concludes ``we have one more reason to entertain
reservations about set theory and that for the time being
no way of rehabilitating this theory is known.'' Solovay
\cite{Solovay2010:PersonalCommunication}, regarding this
reaction as ``rather hysterical,'' nevertheless points out
that, ``[von Neumann's] paper, as a whole, was an important
way station in the evolution of formal set theory. (The
system of \Godel's orange monograph is a direct
descendant.)''

\section{The Ontology of Forcing}
\label{Section:OntologyOfForcing}

A central question in the dispute between the universe view
and the multiverse view is whether there are universes
outside $V$, the universe taken under the universe view to
be the absolute background universe of all sets. A special
case of this question, of course, concerns the status of
forcing extensions of $V$, and this special case appears in
many ways to capture the full debate. On the universe view,
of course, forcing extensions of $V$ are deemed illusory,
for $V$ is already everything, while the multiverse
perspective regards $V$ as a relative concept, referring to
whichever universe is currently under consideration,
without there being any absolute background universe. On
the multiverse view, the use of the symbol $V$ to mean
``the universe'' is something like an introduced constant
that might refer to any of the universes in the multiverse,
and for each of these the corresponding forcing extensions
$V[G]$ are fully real. Thus, I find the following question,
concerning the ontological status of forcing extensions
$V[G]$ of the universe $V$, to be at the heart of the
matter.

\begin{question}
Do forcing extensions of the universe
exist?
\end{question}

More concretely, if $\P\in V$ is a nontrivial forcing
notion, is there a $V$-generic filter $G$ over $\P$? Of
course, all set theorists agree that there can be no such
filter {\it in} $V$. Proponents of the universe view, which
take $V$ to be everything, therefore answer negatively, and
the assertion one sometimes hears
$$\centerline{\it``There are no $V$-generic filters''}$$ is a catechism for the universe view. Such an answer shares a parallel with the assertion
$$\centerline{\it ``There is no square root of $-1$.''}$$ Of course, $\sqrt{-1}$ does not exist in the real field $\R$. One must go to the
field extension, the complex numbers, to find it.
Similarly, one must go to the forcing extension $V[G]$ to
find $G$. Historically, $\sqrt{-1}$ was viewed with
suspicion, and existence deemed imaginary, but useful.
Nevertheless, early mathematicians manipulated expressions
like $2+\sqrt{-5}$ and thereby caught a glimpse into the
richer mathematical world of complex numbers. Similarly,
set theorists now manipulate names and the forcing
relation, thereby catching a glimpse into the forcing
extensions $V[G]$. Eventually, of course, mathematicians
realized how to simulate the complex numbers $a+bi\in\C$
concretely inside the real numbers, representing them as
pairs $(a,b)$ with a peculiar multiplication
$(a,b)\cdot(c,d)=(ac-bd,ad+bc)$. This way, one gains some
access to the complex numbers, or a simulation of them,
from a world having only real numbers, and full acceptance
of complex numbers was on its way. The case of forcing has
some similarities. Although there is no generic filter $G$
inside $V$, there are various ways of simulating the
forcing extension $V[G]$ inside $V$, using the forcing
relation, or using the Boolean-valued structure $V^\B$, or
by using the Naturalist account of forcing. None of these
methods provides a full isomorphic copy of the forcing
extension inside the ground model (as the complex numbers
are simulated in the reals), and indeed they provably
cannot---it is simply too much to ask---but nevertheless
some of the methods come maddeningly close to this. In any
case it is true that with forcing, one has a high degree of
access from the ground model to the forcing extension.

I would like to explain some of the details of this aspect
of the various approaches to forcing, and so let us briefly
review them. The forcing method was introduced by Cohen
\cite{Cohen1963:IndependenceOfCH,Cohen1964:IndependenceOfCHII,Cohen1966:SetTheoryAndTheContinuumHypothesis},
and was used initially to prove the independence of the
axiom of choice and the continuum hypothesis. It was
followed by an explosion of applications that produced an
enormous variety of models of set theory. Forcing revealed
the huge extent of the independence phenomenon, for which
numerous set-theoretic assertions are independent of the
\ZFC\ axioms of set theory and others. With forcing, one
begins with a ground model model $V$ of set theory and a
partial order $\P$ in $V$. Supposing that $G\of\P$ is a
$V$-generic filter, meaning that $G$ contains members from
every dense subset of $\P$ in $V$, one proceeds to build
the forcing extension $V[G]$ by closing under elementary
set-building operations.
$$V\of V[G]$$ In effect, the forcing extension has adjoined the ``ideal'' object $G$ to $V$, in much the same way that one might build a field extension such as
$\Q[\sqrt{2}]$. In particular, every object in $V[G]$ has a
name in $V$ and is constructed algebraically from its name
and $G$. Remarkably, the forcing extension $V[G]$ is always
a model of \ZFC. But it can exhibit different set-theoretic
truths in a way that can be precisely controlled by the
choice of $\P$. The ground model $V$ has a surprising
degree of access to the objects and truths of $V[G]$. The
method has by now produced a staggering collection of
models of set theory.

Some accounts of forcing proceed by defining the {\df
forcing relation} $p\forces\varphi$, which holds whenever
every $V$-generic filter $G$ containing the condition $p$
has $V[G]\satisfies\varphi$. The fundamental facts of
forcing are that: (1) the forcing extension $V[G]$
satisfies \ZFC. (2) Every statement $\varphi$ that holds in
$V[G]$ is forced by some condition $p$ in $G$, and (3) the
forcing relation $p\forces\varphi$ is definable in the
ground model (for fixed $\varphi$, or for $\varphi$ of
fixed complexity).

In order to construct the forcing extension, it appears
that we need to find a suitably generic filter, and how do
we do this? The difficulty of finding a generic filter is,
of course, the crux. One traditional resolution of the
difficulty is to force over a countable transitive ground
model $M$. I shall call this the {\it countable transitive
ground model} method of forcing, and this is also the
initial instance of what I shall later refer to as the {\it
toy model} method of formalization. One starts with such a
countable transitive model $M$. Because it is countable, it
has only countably many dense subsets of $\P$, which we may
enumerate externally as $D_0$, $D_1$, $D_2$, and so on, and
then proceed to pick any condition $p_0\in D_0$, and then
$p_1$ in $D_1$ below $p_0$, and so on. In this way, we
build by diagonalization a descending sequence $p_0\geq
p_1\geq p_2\geq\cdots$, such that $p_n\in D_n$. It follows
that the filter $G$ generated by this sequence is
$M$-generic, and we may proceed to construct the forcing
extension $M[G]$. Thus, with the countable transitive model
approach, one arrives at the forcing extensions in a very
concrete manner. Indeed, the elementary diagram of $M[G]$
is Turing computable from the diagram of $M$ and vice
versa.

There are a number of drawbacks, however, to the countable
transitive ground model approach to forcing. The first
drawback is that it provides an understanding of forcing
over only some models of set theory, whereas other accounts
of forcing allow one to make sense of forcing over any
model of set theory. With the countable transitive model
approach to forcing, for example, the question ``Is
$\varphi$ forceable?'' appears sensible only when asked in
connection with a countable transitive model $M$, and this
is an impoverishment of the method. A second drawback
concerns metamathematical issues surrounding the existence
of countable transitive models of \ZFC: the basic problem
is that we cannot prove that there are any such models,
because by G\"odel's Incompleteness Theorem, if \ZFC\ is
consistent then it cannot prove that there are any models
of \ZFC\ at all. Even if we were to assume $\Con(\ZFC)$,
then we still can't prove that there is a {\it transitive}
model of \ZFC, since the existence of such a model implies
$\Con(\ZFC+\Con(\ZFC))$, and the consistency of this, and
so on transfinitely. In this sense, the existence of a
transitive model of \ZFC\ is something like a very weak
large cardinal axiom, which, although weaker than the
assertion that there is an inaccessible cardinal,
nevertheless transcends \ZFC\ in consistency strength. The
problem recurs with any stronger theory, such as \ZFC\ plus
large cardinals, since no theory (nor the consistency of
the theory) can prove the existence of a transitive model
of that same theory. Thus, it is no help to work in a
background theory with many large cardinals, if one also
wants to do forcing over such models. As a result, this
approach to forcing seems to require one to pay a sort of
tax just to implement the forcing method, starting with a
stronger hypothesis than one ends up with just in order to
carry out the argument. One way to avoid this problem,
which is rarely implemented even though it answers this
objection completely, is simply to drop the transitivity
part of the hypothesis; the transitivity of the ground
model is not actually needed in the construction, provided
one is not squeamish about ill-founded models, for the
countability of $M$ is sufficient to produce the
$M$-generic filter, and with care one can still build the
forcing extension $M[G]$ even when $M$ is ill-founded. This
solution may be rare precisely because set-theorists favor
the well-founded models. More commonly, the defect is
addressed by working not with countable transitive models
of full \ZFC, but rather with countable transitive models
of some finite fragment of \ZFC, whose existence is a
consequence of the Reflection Theorem. Such a fragment
$\ZFC^*$ is usually unspecified, but deemed `sufficiently
large' to carry out all the usual set-theoretic
constructions, and set theorists habitually take $\ZFC^*$
to be essentially a substitute for \ZFC. To prove an
independence result using this method, for a set-theoretic
assertion $\varphi$, one shows that every countable
transitive model $M$ of a sufficient finite fragment
$\ZFC^{**}$ has a forcing extension $M[G]$ satisfying
$\ZFC^*+\varphi$, and concludes as a result that $\ZFC$
cannot prove $\neg\phi$, for if there were such a proof,
then it would use only finitely many of the $\ZFC$ axioms,
and we could therefore put those axioms in a finite
fragment $\ZFC^*$, find a model $M\satisfies\ZFC^{**}$ and
build the corresponding $M[G]\satisfies\ZFC^*+\varphi$, a
contradiction. Such an argument establishes
$\Con(\ZFC)\implies\Con(\ZFC+\varphi)$. Although this move
successfully avoids the extra assumption of the existence
of countable transitive models, it presents drawbacks of
its own, the principal one being that it does not actually
produce a model of $\ZFC+\varphi$, but only shows that this
theory is consistent, leaving one ostensibly to produce the
model by the Henkin construction or some other such method.
Another serious drawback is that this approach to forcing
pushes much of the technique inappropriately into the
meta-theory, for at the step where one says that the proof
of $\neg\varphi$ involves only finitely many axioms and
then appeals to the Reflection Theorem to get the model
$M$, one must have the specific list of axioms used, for we
have the Reflection Theorem only as a theorem scheme, and
not as a single assertion ``For every finite list of
formulas, there is a $V_\alpha$ modeling them...''. Thus,
because the method mixes together the internal notion of
finiteness with that of the meta-theory, the result is that
important parts of the proof take place in the meta-theory
rather than in \ZFC, and this is meta-mathematically
unsatisfying in comparison with the approaches that
formalize forcing entirely as an internal construction.

Other accounts of forcing, in contrast, can be carried out
entirely within \ZFC, allowing one to force over any model
of \ZFC, without requiring one to consider any countable
transitive models or to make any significant moves external
to the model or in the meta-theory. One such method is what
I call the Naturalist account of forcing, which I shall
shortly explain. This method is closely related to (but not
the same as) the Boolean-valued model approach to forcing,
another traditional approach to forcing, which I shall
explain first. The concept of a Boolean-valued model has
nothing especially to do with set theory, and one can have
Boolean-valued partial orders, graphs, groups, rings, and
so on, for any first order theory, including set theory. A
Boolean-valued structure consists of a complete Boolean
algebra $\B$, together with a collection of objects, called
the {\df names}, and an assignment of the atomic formulas
$\sigma=\tau$ and $R(\sigma,\tau)$ to elements of $\B$,
their Boolean values, denoted $\boolval{\sigma=\tau}$ and
$\boolval{R(\sigma,\tau)}$, in such a way that the axioms
of equality are obeyed (a technical requirement). One then
proceeds to extend the Boolean values to all formulas
inductively, by defining
$\boolval{\varphi\wedge\psi}=\boolval{\varphi}\wedge\boolval{\psi}$
and $\boolval{\neg\varphi}=\neg\boolval{\varphi}$, using
the algebraic Boolean operations in $\B$, and
$\boolval{\exists
x\,\varphi(x)}=\bigvee_{\tau}\boolval{\varphi(\tau)}$,
ranging over all names $\tau$, using the fact that $\B$ is
complete. The Boolean-valued structure is said to be {\it
full} if the Boolean value of every existential sentence is
realized by some particular name, so that $\boolval{\exists
x\,\varphi(x)}=\boolval{\varphi(\tau)}$ for some particular
$\tau$. So to specify a Boolean-valued model, one needs a
collection of objects, called the names, and a definition
of the Boolean value on the atomic relations on them. The
Boolean values of more complex assertions follow by
recursion.

To carry this out in the case of set theory, one begins in
a universe $V$, with a complete Boolean algebra $\B$ in
$V$, and defines the class of names inductively, so that
$\tau$ is a $\B$-name if it consists of pairs
$\<\sigma,b>$, where $\sigma$ is a (previously constructed)
$\B$-name and $b\in\B$. The idea is that this name promises
to put the set named by $\sigma$ into the set named by
$\tau$ with Boolean value at least $b$. The atomic values
are defined to implement this idea, while also obeying the
laws of equality, and one therefore eventually arrives at
$V^\B$ as a Boolean-valued structure.  The remarkable fact
is that every axiom of \ZFC\ holds with Boolean value $1$
in $V^\B$. Furthermore, for a fixed formula $\varphi$, the
map $\tau\mapsto\boolval{\varphi(\tau)}$ is definable in
$V$, since $V$ can carry out the recursive definition, and
the structure $V^\B$ is full. Since Boolean values respect
deduction, it follows that for any $b\neq 0$, the
collection of $\varphi$ with $b\leq\boolval{\varphi}$ is
closed under deduction and contains no contradictions
(since these get Boolean value $0$). Thus, to use this
method to prove relative consistency results, one must
merely find a complete Boolean algebra $\B$ in a model of
\ZFC\ such that $\boolval{\varphi}\neq 0$ for the desired
statement $\varphi$. Since the \ZFC\ axioms all get Boolean
value $1$, it follows that $\ZFC+\varphi$ is also
consistent. It is interesting to note that one can give a
complete development of the Boolean-valued model $V^\B$
without ever talking about genericity, generic filters or
dense sets.

Let me turn now to what I call the Naturalist account of
forcing, which seeks to legitimize the actual practice of
forcing, as it is used by set theorists. In any
set-theoretic argument, a set theorist is operating in a
particular universe $V$, conceived as the (current)
universe of all sets, and whenever it is convenient he or
she asserts ``let $G$ be $V$-generic for the forcing notion
$\P$,'' and then proceeds to make an argument in $V[G]$,
while retaining everything that was previously known about
$V$ and basic facts about how $V$ sits inside $V[G]$. This
is the established pattern of how forcing is used in the
literature, and it is proved legitimate by the following
theorem.

\begin{theorem}[Naturalist Account of Forcing] If $V$ is a
(the) universe of set theory and $\P$ is a notion of
forcing, then there is in $V$ a class model of the theory
expressing what it means to be a forcing extension of $V$.
Specifically, in the language with $\in$, constant symbols
for every element of $V$, a predicate for $V$, and constant
symbol $G$, the theory
asserts:\label{Theorem.NaturalistAccount}
\begin{enumerate}
 \item The full elementary diagram of $V$, relativized
     to the predicate for $V$.
 \item The assertion that $V$ is a transitive proper
     class in the (new) universe.
 \item The assertion that $G$ is a $V$-generic
     ultrafilter on $\P$.
 \item The assertion that the (new) universe is $V[G]$,
     and \ZFC\ holds there.
\end{enumerate}
\end{theorem}

This is really a theorem scheme, since $V$ does not have
full uniform access to its own elementary diagram, by
Tarski's theorem on the non-definability of truth. Rather,
the theorem identifies a particular definable class model,
and then asserts as a scheme that it satisfies all the
desired properties. Another way to describe the method is
the following:

\begin{theorem}
For any forcing notion $\P$, there is an
elementary embedding
$$V\precsim\Vbar\of\Vbar[G]$$ of the universe $V$ into a class model $\Vbar$ for which there is a $\Vbar$-generic filter $G\of\bar\P$. In particular, $\Vbar[G]$ is a
forcing extension of $\Vbar$, and the entire extension
$\Vbar[G]$, including the embedding of $V$ into $\Vbar$,
are definable classes in $V$, and $G\in
V$.\label{Theorem.BooleanUltrapower}
\end{theorem}

The point here is that $V$ has full access to the model
$\Vbar[G]$, including the object $G$ and the way that $V$
sits inside $\Vbar$, for they are all definable classes in
$V$. The structure $\Vbar[G]$ is the model asserted to
exist in Theorem \ref{Theorem.NaturalistAccount}, and
$\Vbar$ is the interpretation of the predicate for the
ground model $V$ in that theory. Theorem
\ref{Theorem.BooleanUltrapower} is also a theorem scheme,
since the elementarity of the embedding is made as a
separate assertion of elementarity for each formula. The
models $\Vbar$ and $\Vbar[G]$ are not necessarily
transitive or well-founded, and one should view them as
structures $\<\Vbar,\bar\in>$, $\<\Vbar[G],\bar\in>$ with
their own set membership relation $\bar\in$. Under general
conditions connected with large cardinals, however, there
is a substantial class of cases for which one can arrange
that they are transitive (this is the main theme of
\cite{HamkinsSeabold:BooleanUltrapowers}), and in this
case, the embedding $V\precsim\Vbar$ is, of course, a large
cardinal embedding.

The connection between the Naturalist Account of forcing
and the Boolean-valued model approach to forcing is that
the easiest way to prove Theorems
\ref{Theorem.NaturalistAccount} and
\ref{Theorem.BooleanUltrapower} is by means of
Boolean-valued models and more precisely, the Boolean
ultrapower, a concept going back to Vopenka in the 1960s.
(See \cite{HamkinsSeabold:BooleanUltrapowers} for a full
account of the Boolean ultrapower.) If $\B$ is any complete
Boolean algebra, then one can introduce a predicate $\check
V$ for the ground model into the forcing language, defining
$\boolval{\tau\in\check V}=\bigvee_{x\in
V}\boolval{\tau=\check x}$, where $\check x=\set{\<\check
y,1>\st y\in x}$ is the canonical name for $x$. If $\dot
G=\set{\<\check b,b>\st b\in\B}$ is the canonical name for
the generic filter, then the theory of Theorem
\ref{Theorem.NaturalistAccount} has Boolean value $1$. The
last claim in that theory amounts to the technical
assertion $\boolval{\tau=\val(\check\tau,\dot G)}=1$,
expressing the fact that the structure $V^\B$ knows with
Boolean value $1$ that the set named by $\tau$ is indeed
named by $\tau$ via the generic object.

Let me give a few more details. Any Boolean-valued
structure, whether it is a Boolean-valued graph, a
Boolean-valued partial order, a Boolean-valued group or a
Boolean-valued model of set theory, can be transformed into
a classical 2-valued first-order structure simply by taking
the quotient by an ultrafilter. Specifically, let $U\of\B$
be any ultrafilter on $\B$. There is no need for $U$ to be
generic in any sense, and $U\in V$ is completely fine.
Define an equivalence relation on the names by
$\sigma\equiv_U\tau\iff\boolval{\sigma=\tau}\in U$. A
subtle and sometimes misunderstood point is that this is
not necessarily the same as $\val(\sigma,U)=\val(\tau,U)$,
when $U$ is not $V$-generic. Nevertheless, the relation
$\sigma\in_U\tau\iff\boolval{\sigma\in\tau}$ is
well-defined on the equivalence classes, and one can form
the quotient structure $V^\B/U$, using Scott's trick
concerning reduced equivalence classes. The quotient
$V^\B/U$ is now a classical $2$-valued structure, and one
proves the \L os theorem that
$V^\B/U\satisfies\varphi(\tau)\iff\boolval{\varphi(\tau)}\in
U$. The collection $\Vbar$ of (equivalence classes of)
names $\tau$ with $\boolval{\tau\in\check V}\in U$ serves
as the ground model of $V^\B/U$, and one proves that
$V^\B/U$ is precisely $\Vbar[G]$, where $G$ is the
equivalence class $[\dot G]_U$. That is, the object $G$ is
merely the equivalence class of the canonical name for the
generic object, and this is a set that exists already in
$V$. The entire structure $\Vbar[G]$ therefore exists as a
class in $V$, if the ultrafilter $U$ is in $V$. The map
$x\mapsto[\check x]_U$ is an elementary embedding from $V$
to $\Vbar$, the map mentioned in Theorem
\ref{Theorem.BooleanUltrapower}. This map is known as the
Boolean ultrapower map, and it provides a way to generalize
the ultrapower concept from ultrapowers by ultrafilters on
power sets to ultrapowers by ultrafilters on any complete
Boolean algebra. One can equivalently formulate the Boolean
ultrapower map as the direct limit of the system of
classical ultrapowers, obtained via the ultrafilters
generated by $U$ on the power sets of all the various
maximal antichains in $\B$, ordering the antichains under
refinement. These induced ultrapower embeddings form a
large commutative diagram, of which the Boolean ultrapower
is the direct limit.

Part of the attraction of the Naturalist Account of forcing
as developed in Theorem \ref{Theorem.NaturalistAccount} is
that one may invoke the theorem without paying attention to
the manner in which it was proved. In such an application,
one exists inside a universe $V$, currently thought of as
the universe of all sets, and then, invoking Theorem
\ref{Theorem.NaturalistAccount} via the assertion
$$\centerline{``Let $G\of\B$ be $V$-generic. Argue in $V[G]$...''}$$
one adopts the new theory of Theorem
\ref{Theorem.NaturalistAccount}. The theory explicitly
stated in Theorem \ref{Theorem.NaturalistAccount} allows
one to keep all the previous knowledge about $V$,
relativized to a predicate for $V$, but adopt the new (now
current) universe $V[G]$, a forcing extension of $V$. Thus,
although the proof did not provide an actual $V$-generic
filter, the effect of the new theory is entirely as if it
had. This method of application, therefore, implements in
effect the content of the multiverse view. That is, whether
or not the forcing extensions of $V$ actually exist, we are
able to behave via the naturalist account of forcing
entirely as if they do. In any set-theoretic context,
whatever the current set-theoretic background universe $V$,
one may at any time use forcing to jump to a universe
$V[G]$ having a $V$-generic filter $G$, and this jump
corresponds to an invocation of Theorem
\ref{Theorem.NaturalistAccount}.

Of course, one might on the universe view simply use the
naturalist account of forcing as the means of explaining
the illusion:  the forcing extensions don't really exist,
but the naturalist account merely makes it seem as though
they do. The multiverse view, however, takes this use of
forcing at face value, as evidence that there actually are
$V$-generic filters and the corresponding universes $V[G]$
to which they give rise, existing outside the universe.
This is a claim that we cannot prove within set theory, of
course, but the philosophical position makes sense of our
experience---in a way that the universe view does
not---simply by filling in the gaps, by positing as a
philosophical claim the actual existence of the generic
objects which forcing comes so close to grasping, without
actually grasping. With forcing, we seem to have discovered
the existence of other mathematical universes, outside our
own universe, and the multiverse view asserts that yes,
indeed, this is the case. We have access to these
extensions via names and the forcing relation, even though
this access is imperfect. Like Galileo, peering through his
telescope at the moons of Jupiter and inferring the
existence of other worlds, catching a glimpse of what it
would be like to live on them, set theorists have seen via
forcing that divergent concepts of set lead to new
set-theoretic worlds, extending our previous universe, and
many are now busy studying what it would be like to live in
them.

An interesting change in the use of forcing has occurred in
the history of set theory. In the earlier days of forcing,
the method was used principally to prove independence
results, with theorems usually having the form
$\Con(\ZFC+\varphi)\implies\Con(\ZFC+\psi)$, proved by
starting with a model of $\varphi$ and providing $\psi$ in
a forcing extension. Contemporary work would state the
theorem as: {\sl If $\varphi$, then there is a forcing
extension with $\psi$.} Furthermore, one would go on to
explain what kind of forcing was involved, whether it was
cardinal-preserving or c.c.c. or proper and so on. Such a
description of the situation retains important information
connecting the two models, and by emphasizing the relation
between the two models, this manner of presentation
conforms with the multiverse view.

\section{The analogy between set theory and geometry}
\label{Section:AnalogyWithGeometry}

There is a very strong analogy between the multiverse view
in set theory and the most commonly held views about the
nature of geometry. For two thousand years, mathematicians
studied geometry, proving theorems about and making
constructions in what seemed to be the unique background
geometrical universe. In the late nineteenth century,
however, geometers were shocked to discover non-Euclidean
geometries. At first, these alternative geometries were
presented merely as simulations within Euclidean geometry,
as a kind of playful or temporary re-interpretation of the
basic geometric concepts. For example, by temporarily
regarding `line' to mean a great circle on the unit sphere,
one arrives at spherical geometry, where all lines
intersect; by next regarding `line' to mean a circle
perpendicular to the unit circle, one arrives at one of the
hyperbolic geometries, where there are many parallels to a
given line through a given point. At first, these
alternative geometries were considered as curiosities,
useful perhaps for independence results, for with them one
can prove that the parallel postulate is not provable from
the other axioms. In time, however, geometers gained
experience in the alternative geometries, developing
intuitions about what it is like to live in them, and
gradually they accepted the alternatives as geometrically
meaningful. Today, geometers have a deep understanding of
the alternative geometries, which are regarded as fully
real and geometrical.

The situation with set theory is the same. The initial
concept of set put forth by Cantor and developed in the
early days of set theory seemed to be about a unique
concept of set, with set-theoretic arguments and
constructions seeming to take place in a unique background
set-theoretic universe. Beginning with G\"odel's
constructible universe $L$ and particularly with the rise
of forcing, however, alternative set theoretic universes
became known, and today set theory is saturated with them.
Like the initial reactions to non-Euclidean geometry, the
universe view regards these alternative universes as not
fully real, while granting their usefulness for proving
independence results. Meanwhile, set theorists continued,
like the geometers a century ago, to gain experience living
in the alternative set-theoretic worlds, and the multiverse
view now makes the same step in set theory that geometers
ultimately made long ago, namely, to accept the alternative
worlds as fully real.

The analogy between geometry and set theory extends to the
mathematical details about how one reasons about the
alternative geometric and set-theoretic worlds. Geometers
study the alternative geometries in several ways. First,
they study a geometry by means of a simulation of it within
Euclidean space, as in the presentations of spherical and
hyperbolic geometry I mentioned above. This is studying a
geometrical universe from the perspective of another
geometrical universe that has some access to it or to a
simulation of it. Second, they can in a sense jump inside
the alternative geometry, for example by adopting
particular negations of the parallel postulate and
reasoning totally within that new geometrical system.
Finally, third, in a sophisticated contemporary
understanding, they can reason abstractly by using the
group of isometries that defines a particular geometry. The
case of forcing offers these same three modes of reasoning.
First, set theorists may reason about a forcing extension
from the perspective of the ground model via names and the
forcing relation. Second, they may reason about the forcing
extension by jumping into it and reasoning as though they
were living in that extension. Finally, third, they may
reason about the forcing extension abstractly, by examining
the Boolean algebra and its automorphism group, considering
homogeneity properties, in order to make conclusions about
the forcing extension.

A stubborn geometer might insist---like an
exotic-travelogue writer who never actually ventures west
of seventh avenue---that only Euclidean geometry is real
and that all the various non-Euclidean geometries are
merely curious simulations within it. Such a position is
self-consistent, although stifling, for it appears to miss
out on the geometrical insights that can arise from the
other modes of reasoning. Similarly, a set theorist with
the universe view can insist on an absolute background
universe $V$, regarding all forcing extensions and other
models as curious complex simulations within it. (I have
personally witnessed the necessary contortions for class
forcing.) Such a perspective may be entirely
self-consistent, and I am not arguing that the universe
view is incoherent, but rather, my point is that if one
regards all outer models of the universe as merely
simulated inside it via complex formalisms, one may miss
out on insights that could arise from the simpler
philosophical attitude taking them as fully real.

The history of mathematics provides numerous examples where
initially puzzling imaginary objects become accepted as
real. Irrational numbers, such as $\sqrt{2}$, became
accepted; zero became a number; negative numbers and then
imaginary and complex numbers were accepted; and then
non-Euclidean geometries. Now is the time for $V$-generic
filters.

\section{Multiverse response to the categoricity arguments}
\label{Section:Categoricity}

There is a rich tradition of categoricity arguments in set
theory, going back to Peano's
\cite{Peano1889:PrinciplesOfArithmeticPresentedByNewMethod}
\cite[p.
83--97]{vanHeijenoort1967:FromFregeToGodelSourcebookInMathLogic}
proof that his second-order axioms characterize the unique
structure of the natural numbers and to Zermelo's
\cite{Zermelo1930:UberGrenzzahlenUndMengenBereiche}
second-order categoricity proof for set theory,
establishing the possible universes as $V_\kappa$ where
$\kappa$ is inaccessible.\footnote{In what appears to be an
interesting case of convergent evolution in the foundations
of mathematics, this latter universe concept coincides
almost completely with the concept of Grothendieck
universe, now pervasively used in category theory
\cite{Kromer2001:TarskisAxiomOfInaccessiblesAndGrothendieckUniverses}.
The only difference is that the category theorists also
view $V_\omega$ as a Grothendieck universe, which amounts
to considering $\aleph_0$ as an incipient inaccessible
cardinal. Surely the rise of Grothedieck universes in
category theory shares strong affinities with the
multiverse view in set theory, although most set theorists
find Grothendieck universes clumsy in comparison with the
more flexible concept of a (transitive) model of set
theory; nevertheless, the category theorists will point to
the multiverse concepts present in the theory of topuses as
more general still (see
\cite{Blass1984:InteractionBetweenCategoryTheoryAndSetTheory}).}
Continuing the categoricity tradition, Martin
\cite{Martin2001:MultipleUniversesOfSetsAndIndeterminateTruthValues}
argues, to explain it very briefly, that the set-theoretic
universe is unique, because any two set-theoretic concepts
$V$ and $V'$ can be compared level-by-level through the
ordinals, and at each stage, if they agree on
$V_\alpha=V'_\alpha$ and each is claiming to be all the
sets, then they will agree on $V_{\alpha+1}=V'_{\alpha+1}$,
and so ultimately $V=V'$. The categoricity arguments, of
course, tend to support the universe view.

The multiversist objects to Martin's presumption that we
are able to compare the two set concepts in a coherent way.
Which set concept are we using when undertaking the
comparison? Martin's argument employs a background concept
of `property,' which amounts to a common set-theoretic
context in which we may simultaneously refer to both set
concepts when performing the inductive comparison. Perhaps
one would want to use either of the set concepts as the
background context for the comparison, but it seems
unwarranted to presume that either of the set concepts is
able to refer to the other internally, and the ability to
make external set (or property) concepts internal is the
key to the success of the induction.

If we make explicit the role of the background
set-theoretic context, then the argument appears to reduce
to the claim that {\it within any fixed set-theoretic
background concept}, any set concept that has all the sets
agrees with that background concept; and hence any two of
them agree with each other. But such a claim seems far from
categoricity, should one entertain the idea that there can
be different incompatible set-theoretic backgrounds.

Another more specific issue with Martin's comparison is
that it requires that the two set concepts agree on the
notion of ordinal sufficiently that one can carry out the
transfinite recursive comparison. To be sure, Martin does
explicitly assume in his argument that the concept of
natural number is sharp in his sense, and he adopts a sharp
account of the concept of wellordering, while admitting
that ``it is of course possible to have doubts about the
sharpness of the concept of wellordering.'' Lacking a fixed
background concept of set or of well-order, it seems that
the comparison may become eventually incoherent. For
example, perhaps each set concept individually provides an
internally coherent account of set and of well-order, but
in any common background, both are revealed as inadequate.

There seems little reason why two different concepts of set
need to agree even on the concept of the natural numbers.
Although we conventionally describe the natural numbers as
$1$, $2$, $3$, $\ldots$, and so on, why are we so confident
that this ellipses is meaningful as an absolute
characterization? Peano's categoricity proof is a
second-order proof that is sensible only in the context of
a fixed concept of subsets of $\mathbb{N}$, and so this
ellipses carries the baggage of a set-theoretic ontology.
Our initial confidence that our ``and so on'' describes a
unique structure of natural numbers should be tempered by
the comparatively opaque task of grasping all sets of
natural numbers, whose existence and nature supports and is
required for the categoricity proof. On the multiverse
view, the  possibility that differing set concepts may lead
to different and perhaps incomparable concepts of natural
number suggests that Martin's comparison process may
eventually become incoherent.

To illustrate, consider the situation of Peano's
categoricity result for the natural numbers. Set theory
surely provides a natural context in which to carry out
Peano's second-order argument that all models of
second-order arithmetic are canonically isomorphic. Indeed,
one may prove in \ZFC\ that there is a unique
(second-order) inductive structure of the natural numbers.
So once we fix a background concept of set satisfying (a
small fragment of) \ZFC, we achieve Peano's categoricity
result for the natural numbers. Nevertheless, we also know
quite well that different models of set theory can provide
different and quite incompatible background concepts of
set, and different models of \ZFC\ can have quite different
incomparable versions of their respective standard $\N$s.
Even though we prove in \ZFC\ that $\N$ is unique, the
situation is that not all models of \ZFC\ have the same
$\N$ or even the same arithmetic truths. Victoria Gitman
and I have defined that a model of arithmetic is {\df a}
standard model of arithmetic, as opposed to {\it the}
standard model of arithmetic, if it is the $\N$ of a model
of \ZFC, and Ali Enayat has characterized the nonstandard
instances as precisely the computably-saturated models of
the arithmetic consequences of \ZFC, and indeed, these are
axiomatized by \PA\ plus the assertions
$\varphi\to\Con(\ZFC_n+\varphi^{\N})$, where $\ZFC_n$ is
the $\Sigma_n$ fragment of $\ZFC$.

The point is that a second-order categoricity argument,
even just for the natural numbers, requires one to operate
in a context with a background concept of set. And so
although it may seem that saying ``$1$, $2$, $3$, $\ldots$
and so on,'' has to do only with a highly absolute concept
of finite number, the fact that the structure of the finite
numbers is uniquely determined depends on our much murkier
understanding of which subsets of the natural numbers
exist. So why are mathematicians so confident that there is
an absolute concept of finite natural number, independent
of any set theoretic concerns, when all of our categoricity
arguments are explicitly set-theoretic and require one to
commit to a background concept of set? My long-term
expectation is that technical developments will eventually
arise that provide a forcing analogue for arithmetic,
allowing us to modify diverse models of arithmetic in a
fundamental and flexible way, just as we now modify models
of set theory by forcing, and this development will
challenge our confidence in the uniqueness of the natural
number structure, just as set-theoretic forcing has
challenged our confidence in a unique absolute
set-theoretic universe.

\section{The multiverse provides a context for universe adjudication}
\label{Section:Adjudication}

From the multiverse perspective, the advocates of the
universe view want in effect to fix a particular set
concept, with associated particular universe $V$, and
declare it as the absolute background concept of set. In
this way, the multiverse view can simulate the universe
view by restricting attention to the lower cone in the
multiverse consisting of this particular universe $V$ and
the universes below it. Inside this cone, the universe $V$
provides in effect absolute background notions of
countability, well-foundedness, and so on, and of course
there are no nontrivial $V$-generic filters in this
restricted multiverse. Any multiverse set-theorist can
pretend to be a universe set-theorist simply by jumping
into a specific $V$ and temporarily forgetting about the
worlds outside $V$.

Meanwhile, there is understandably a level of disagreement
within the universe community about precisely which
universe $V$ to fix, about the features and truths of the
{\it real} universe. The various arguments about the final
answers to set-theoretic questions---is the continuum
hypothesis {\it really} true or not?---amount to questions
about which or what kind of $V$ we shall fix as the
absolute background. The point here is that this is a
debate that naturally takes place within the multiverse
arena. It is explicitly a multiverse task to compare
differing set concepts or universes or to entertain the
possibility that there are fundamentally differing but
legitimate concepts of set. It is only in a multiverse
context that we may sensibly compare competing proposals
for the unique absolute background universe. In this sense,
the multiverse view provides a natural forum in which to
adjudicate differences of opinion arising within the
universe view.

\section{Case study: multiverse view on the continuum hypothesis}
\label{Section:CaseStudyCH}

Let me discuss, as a kind of case study, how the multiverse
and universe views treat one of the most important problems
in set theory, the continuum hypothesis (\CH). This is the
famous assertion that every set of reals is either
countable or equinumerous with $\R$, and it was a major
open question from the time of Cantor, appearing at the top
of Hilbert's well-known list of open problems in 1900. The
continuum hypothesis is now known to be neither provable
nor refutable from the usual \ZFC\ axioms of set theory, if
these axioms are consistent. Specifically, G\"odel proved
that $\ZFC+\CH$ holds in the constructible universe $L$ of
any model of $\ZFC$, and so $\CH$ is not refutable in
$\ZFC$. In contrast, Cohen proved that $L$ has a forcing
extension $L[G]$ satisfying $\ZFC+\neg\CH$, and so \CH\ is
not provable in \ZFC. The generic filter directly adds any
number of new real numbers, so that there could be
$\aleph_2$ of them or more in $L[G]$, violating \CH.

Going well beyond mere independence, however, it turns out
that both \CH\ and $\neg\CH$ are forceable over any model
of set theory.

\begin{theorem} The universe $V$ has forcing extensions
\begin{enumerate}
 \item $V[G]$, collapsing no cardinals, such that
     $V[G]\satisfies\neg\CH$.
 \item $V[H]$, adding no new reals, such that
     $V[H]\satisfies\CH$.
\end{enumerate}
\end{theorem}

In this sense, every model of set theory is very close to
models with the opposite answer to \CH. Since the $\CH$ and
$\neg\CH$ are easily forceable, the continuum hypothesis is
something like a lightswitch, which can be turned on and
off by moving to ever larger forcing extensions. In each
case the forcing is relatively mild, with the new
universes, for example, having all the same large cardinals
as the original universe. After decades of experience and
study, set-theorists now have a profound understanding of
how to achieve the continuum hypothesis or its negation in
diverse models of set theory---forcing it or its negation
in innumerable ways, while simultaneously controlling other
set-theoretic properties---and have therefore come to a
deep knowledge of the extent of the continuum hypothesis
and its negation in the multiverse.

On the multiverse view, consequently, the continuum
hypothesis is a settled question; it is incorrect to
describe the \CH\ as an open problem. The answer to \CH\
consists of the expansive, detailed knowledge set theorists
have gained about the extent to which it holds and fails in
the multiverse, about how to achieve it or its negation in
combination with other diverse set-theoretic properties. Of
course, there are and will always remain questions about
whether one can achieve \CH\ or its negation with this or
that hypothesis, but the point is that the most important
and essential facts about \CH\ are deeply understood, and
these facts constitute the answer to the \CH\ question.

To buttress this claim, let me offer a brief philosophical
argument that the \CH\ can no longer be settled in the
manner that set theorists formerly hoped that it might be.
Set theorists traditionally hoped to settle \CH\ according
to the following template, which I shall now refer to as
the {\em dream solution} template for \CH:

\begin{quote}
\medskip
{\bf Step 1.} Produce a set-theoretic assertion $\Phi$
expressing a natural `obviously true' set-theoretic
principle.

\medskip\noindent{\bf Step 2.} Prove that $\Phi$ determines \CH.

That is, prove that $\Phi\implies\CH$, or prove that
$\Phi\implies\neg\CH$.
\end{quote}

\noindent In step 1, the assertion $\Phi$ should be
obviously true in the same sense that many set-theorists
find the axiom of choice and other set-theoretic axioms,
such as the axiom of replacement, to be obviously true,
namely, the statement should express a set-theoretic
principle that we agree should be true in the intended
intended interpretation, the pre-reflective set theory of
our imagination. Succeeding in the dream solution would
settle the \CH, of course, because everyone would accept
$\Phi$ and its consequences, and these consequences would
include either $\CH$ or $\neg\CH$.

I claim now that this dream solution has become impossible.
It will never be realized. The reason has to do with our
rich experience in set-theoretic worlds having \CH\ and
others having $\neg\CH$. Our situation, after all, is not
merely that \CH\ is formally independent and we have no
additional knowledge about whether it is true or not.
Rather, we have an informed, deep understanding of the \CH\
and $\neg\CH$ worlds and of how to build them from each
other. Set theorists today grew up in these worlds, and
they have flicked the \CH\ light switch many times in order
to achieve various set-theoretic effects. Consequently, if
you were to present a principle $\Phi$ and prove that it
implies $\neg\CH$, say, then we can no longer see $\Phi$ as
obviously true, since to do so would negate our experiences
in the set-theoretic worlds having $\CH$. Similarly, if
$\Phi$ were proved to imply $\CH$, then we would not accept
it as obviously true, since this would negate our
experiences in the worlds having $\neg\CH$. The situation
would be like having a purported `obviously true' principle
that implied that midtown Manhattan doesn't exist. But I
know it exists; I live there. Please come visit! Similarly,
both the \CH\ and $\neg\CH$ worlds in which we have lived
and worked seem perfectly legitimate and fully
set-theoretic to us, and because of this, either
implication in step 2 immediately casts doubt to us on the
naturality of $\Phi$. Once we learn that a principle
fulfills step 2, we can no longer accept it as fulfilling
step 1, even if previously we might have thought it did. My
predicted response to any attempt to carry out the dream
solution is that claims of mathematical naturality (in step
1) will be met by objections arising from deep mathematical
experience of the contrary situations (in step 2).

Allow me to give a few examples that illustrate this
prediction. The first is the response to Freiling's
\cite{Freiling1986:AxiomsOfSymmetry:ThrowingDartsAtTheRealLine}
delightful axiom of symmetry, which he presents as part of
``a simple philosophical `proof' of the negation of
Cantor's continuum hypothesis.'' Freiling spends the first
several pages of his article ``subject[ing] the continuum
to certain thought experiments involving random darts,''
before ultimately landing at his axiom of symmetry, which
he presents as an ``intuitively clear axiom.'' The axiom
itself asserts that for any function $f$ mapping reals to
countable sets of reals, there must reals $x$ and $y$ such
that $y\notin f(x)$ and $x\notin f(y)$. Freiling argues at
some length for the natural appeal of his axiom, asking us
to imagine throwing two darts at a dart board, and having
thrown the first dart, which lands at some position $x$,
observing then that because $f(x)$ is a countable set, we
should expect almost surely that the second dart will land
at a point $y$ not in $f(x)$. And since the order in which
we consider the darts shouldn't seem to matter, we conclude
by symmetry that almost surely $x\notin f(y)$ as well. We
therefore have natural reason, he argues, to expect not
only that there is a pair $(x,y)$ with the desired
property, but that intuitively almost all pairs have the
desired property. Indeed, he asserts that, ``actually [the
axiom], being weaker than our intuition, does not say that
the two darts have to do anything. All it claims is that
what heuristically will happen every time, can happen.'' In
this way, Freiling appears directly to be carrying out step
1 in the dream solution template. And sure enough, he
proceeds next to carry out step 2, by proving in \ZFC\ that
the axiom of symmetry is exactly equivalent to $\neg\CH$.
The forward implication is easy, for if there is a
well-ordering of $\R$ in order type $\omega_1$, then we may
consider the function $f$ mapping every real $x$ to the
initial segment of the order up to $x$, which is a
countable set, and observe by linearity that for any pair
of numbers either $x$ precedes $y$ or conversely, and so
either $x\in f(y)$ or $y\in f(x)$, contrary to the axiom of
symmetry, which establishes $\AS\implies\neg\CH$.
Conversely, if $\CH$ fails, then for any choice of
$\omega_1$ many distinct $x_\alpha$, for $\alpha<\omega_1$,
there must by $\neg\CH$ be a real $y\notin\Union_\alpha
f(x_\alpha)$, since this union has size $\omega_1$, but
$f(y)$ can only contain countably many $x_\alpha$, and so
some late enough $x_\alpha$ and $y$ are as desired. In
summary, $\AS\iff\neg\CH$, and Freiling has exactly carried
out the dream solution template for \CH.

But was his argument received as a solution of \CH? No.
Many mathematicians objected that Freiling was implicitly
assuming for a given function $f$ that various sets were
measurable, including importantly the set $\set{(x,y)\st
y\in f(x)}$. Freiling clearly anticipated this objection,
making the counterargument in his paper that he was
justifying his axioms prior to any mathematical development
of measure, on the same philosophical or pre-reflective
ideas that are used to justify our mathematical
requirements for measure in the first place. Thus, he
argues, we would seem to have as much intuitive support for
the axiom of symmetry directly as we have for the idea that
our measure should be countably additive, for example, or
that it should satisfy the other basic properties we demand
of it.

My point is not to defend Freiling specifically, but rather
to observe that mathematicians objected to Freiling's
argument largely from a perspective of deep experience and
familiarity with non-measurable sets and functions,
including extreme violations of the Fubini property, and
for mathematicians with this experience and familiarity,
the pre-reflective arguments simply fell flat. We are
skeptical of any intuitive or naive use of measure
precisely because we know so much now about the various
mathematical pitfalls that can arise, about how complicated
and badly-behaved functions and sets of reals can be in
terms of their measure-theoretic properties. We know that
any naive account of measure will have a fundamental
problem dealing with subsets of the plane all of whose
horizontal sections are countable and all of whose vertical
sections are co-countable, for example, precisely because
the sets looks very small from one direction and very large
from another direction, while we expect that rotating a set
should not change its size. Indeed, one might imagine a
variant of the Freiling argument, proceeding like this:
intuitively all sets are measurable, and also rotating sets
in the plane preserves measure; but if \CH\ holds, then
there are sets, such as the graph of a well-ordering of
$\R$ in order type $\omega_1$, that look very small from
one direction and very large from another. Hence, \CH\
fails. Nevertheless, our response to this argument is
similar to our response to the original Freiling argument,
namely, because of our detailed experience with badly
behaved sets and functions with respect to their measure
properties, we simply cannot accept the naive expectations
of the modified argument, and in the original case, we
simply are not convinced by Freiling's argument that \AS\
is intuitively true, even if he is using the same
intuitions that guided us to the basic principles of
measure in the first place. In an extreme instance of this,
inverting Freiling's argument, set theorists sometimes
reject the principle as a fundamental axiom precisely
because of the counterexamples to it that one can produce
under \CH.

In the end, Freiling's argument is not generally accepted
as a solution to \CH, and his axiom instead is most often
described as an intriguing equivalent formulation of
$\neg\CH$, which is interesting to consider and which
reveals at its heart the need for us to take care with
issues of non-measurability. In this way, Freiling's
argument is turned on its head, as a warning about the
error that may arise from a naive treatment of measure
concepts. In summary, the episode exhibits my predicted
response for any attempted answer following the dream
solution template: rejection of the new axiom from a
perspective of deep mathematical experience with the
contrary.

Let me turn now to a second example. Consider the
set-theoretic principle that I shall call the {\df powerset
size axiom} {\rm PSA}, which asserts in brief that strictly
larger sets have strictly more subsets:
$$\forall x,y\qquad |x|\lt |y|\implies |P(x)|\lt |P(y)|.$$
Although set-theorists understand the situation of this
axiom very well, as I will shortly explain, allow me first
to discuss how it is received in non-logic and
non-set-theoretic mathematical circles, which is: extremely
well! An enormous number of mathematicians, including many
very good ones, view the axiom as extremely natural or even
obviously true in the same way that various formulations of
the axiom of choice or the other basic principles of set
theory are obviously true. The principle, for example, is
currently the top-rated answer among dozens to a popular
mathoverflow question seeking examples of
reasonable-sounding statements that are nevertheless
independent of the axioms of set theory
\cite{HamkinsMO6594:AnswerToWhatAreSomeReasonable-soundingStatementsIndepOfZFC?}.\footnote{There
are, in addition, at least three other mathoverflow
questions posted by mathematicians simply asking naively
whether the PSA is true, or how to prove it, or indeed
whether it is true that the assertion is not provable,
incredible as that might seem.} In my experience, a brief
talk with mathematicians at your favorite math tea has a
good chance to turn up additional instances of
mathematicians who find the assertion to express a basic
fact about sets.

Meanwhile, set theorists almost never agree with this
assessment, for they know that one can achieve all kinds of
crazy patterns for the continuum function $\kappa\mapsto
2^\kappa$ via Easton's theorem, and even Cohen's original
model of $\ZFC+\neg\CH$ had $2^\omega=2^{\omega_1}$, the
assertion known as Luzin's hypothesis
\cite{Luzin1935:SurLesEnssemblesAnalytiquesNuls}, which had
been proposed as an alternative to the continuum
hypothesis. Furthermore, Martin's axiom implies
$2^\omega=2^\kappa$ for all $\kappa<2^\omega$, which can
mean additional violations of $\text{PSA}$ when \CH\ fails
badly. So not only do set-theorists know that $\text{PSA}$
can fail, but also they know that $\text{PSA}$ must fail in
models of the axioms, such as the proper forcing axiom
\PFA\ or Martin's maximum \MM, that are often favored
particularly by set-theorists with a universe view.

Thus, the situation is that a set-theoretic principle that
many mathematicians find to be obviously true and which
surely expresses an intuitively clear pre-reflective
principle about the concept of size, and which furthermore
is known by set-theorists to be perfectly safe in the sense
that it is relatively consistent with the other axioms of
\ZFC\ and in fact a consequence of the generalized
continuum hypothesis, is nevertheless almost universally
rejected by set-theorists when it is proposed as a
fundamental axiom. This rejection follows my predicted
pattern: we simply know too much about the various ways
that the principle can be violated and have too much
experience working in models of set theory where the
principle fails to accept it as a fundamental axiom.

Imagine briefly that the history of set theory had
proceeded differently; imagine, for example, that the
powerset size axiom had been considered at the very
beginning of set theory---perhaps it was used in the proof
of a critical theorem settling a big open question at the
time---and was subsequently added to the standard list of
axioms. In this case, perhaps we would now look upon models
of $\neg\text{PSA}$ as strange in some fundamental way,
violating a basic intuitive principle of sets concerning
the relative sizes of power sets; perhaps our reaction to
these models would be like the current reaction some
mathematicians (not all) have to models of $\ZF+\neg\AC$ or
to models of Aczel's anti-foundation axiom \AFA, namely,
the view that the models may be interesting mathematically
and useful for a purpose, but ultimately they violate a
basic principle of sets.

There have been other more sophisticated attempts to settle
\CH, which do not rely on the dream solution template.
Woodin, for example, has advanced an argument for $\neg\CH$
based on $\Omega$-logic, appealing to desirable structural
properties of this logic, and more recently has proposed to
settle \CH\ positively in light of his proposed
construction of {\df ultimate-$L$}. In both cases, these
arguments attempt to settle \CH\ by pointing out (in both
cases under conjectured technical assumptions) that in
order for the absolute background universe to be regular or
desirable in a certain very precise, technical way, the
favored solution of \CH\ must hold. To my way of thinking,
however, a similar objection still applies to these
arguments. What the multiversist desires in such a line of
reasoning is an explicit explanation of how our experience
in the \CH\ worlds or in the $\neg\CH$ worlds was somehow
illusory, as it seems it must have been, if the argument is
to succeed. Since we have an informed, deep understanding
of how it could be that \CH\ holds or fails, even in worlds
very close to any given world, it will be difficult to
regard these worlds as imaginary. To sloganize the dispute,
the universists present a shining city on a hill; and
indeed it is lovely, but the multiversists are widely
traveled and know it is not the only one.

\section{Case study: multiverse view on the axiom of
constructibility}\label{Section:CaseStudyV=L}

Let me turn now, for another small case study, to the axiom
of constructibility $V=L$ and how it is treated under the
universe and multiverse positions. The constructible
universe $L$, introduced by G\"odel, is the result of the
familiar proper class length transfinite recursive
procedure $$L=\Union_{\alpha\in\ORD}L_\alpha,$$ where each
$L_\alpha$ consists of the sets explicitly definable from
parameters over the previous levels. The objects in the
constructible universe would seem to be among the minimal
collection of objects that one would be committed to, once
one is living in a universe with those particular ordinals.
The remarkable fact, proved by \Godel, is that because of
its uniform structure, $L$ satisfies all of the \ZF\ axioms
of set theory, as well as the axiom of choice and the
generalized continuum hypothesis; indeed, the fine
structure and condensation properties of $L$ lead to
numerous structural features, such as the $\Diamond$
principle and many others, and as a result the $V=L$
hypothesis settles many questions in set theory. For many
of these questions, this was how set-theorists first knew
that they were relatively consistent, and in this sense,
the axiom $V=L$ has very high explanatory force.

Nevertheless, most set theorists reject $V=L$ as a
fundamental axiom. It is viewed as limiting: why should we
have only constructible sets? Furthermore, although the
$V=L$ hypothesis settles many set-theoretic questions, it
turns out too often to settle them in the ``wrong'' way, in
the way of counterexamples and obstacles, rather than in a
progressively unifying theory. One might compare the
cumbersome descriptive set theory of $L$ to the smooth
uniformity that is provable from projective determinacy. A
significantly related observation is that $L$ can have none
of the larger large cardinals.

Several set theorists have emphasized that even under
$V\neq L$, we may still analyze what life would be like in
$L$ by means of the the translation
$\varphi\mapsto\varphi^L$, that is, by asking in $V$
whether a given statement $\varphi$ is true in the $L$ of
$V$. Steel makes the point as follows in slides from his
talk
\cite{Steel2004:SlidesGenericAbsolutenessAndTheContinuumProblem}
as follows (see my criticism below):

\begin{quotation}$V=L$ is restrictive, in that adopting
it limits the interpretative power of our language.
\begin{enumerate}
\item[1.] The language of set theory as used by the
    $V=L$ believer can be translated into the language
    of set theory as used by the ``$\exists$ measurable
    cardinal'' believer: $$\varphi\mapsto \varphi^L.$$
\item[2.] There is no translation in the other
    direction.
\item[3.] Proving $\varphi$ in $\ZFC+ V=L$
is equivalent to proving $\varphi^L$ in \ZFC. Thus
adding $V=L$ settles no
    questions not already settled by \ZFC. It just prevents
    us from asking as many questions!
\end{enumerate}
\end{quotation}

\noindent Despite the near universal rejection of $V=L$ as
a fundamental axiom, the universe $L$ is intensively
studied, along with many other $L$-like universes
accommodating larger large cardinals. This practice exactly
illustrates the multiverse perspective, for set-theorists
study $L$ as one of many many possible set-theoretic
universes, and our knowledge of $L$ and the other canonical
inner models is particularly interesting when it sheds
light on $V$. At the same time, the multiverse perspective
on $L$ takes note of several observations against the view
of $L$ as fundamentally limiting, for in several surprising
ways, the constructible universe $L$ is not as impoverished
as it is commonly described to be.

To begin, of course, there is not just one $L$. Each
universe of set theory has its own constructible universe,
and these can vary widely: different ordinals, different
large cardinals, different reals and even different
arithmetic truths. On the multiverse view, one should
understand talk of ``the'' constructible universe as
referring to a rough grouping of very different models, all
satisfying $V=L$.

One eye-opening observation is that $V$ and $L$ have
transitive models of exactly the same (constructible)
theories. For example, there is a transitive model of the
theory $\ZFC+$``there is a proper class of supercompact
cardinals'' if and only if there is such a transitive model
inside $L$. The reason is that the assertion that a given
theory $T$ has a transitive model has complexity
$\Sigma^1_2(T)$, in the form ``there is a real coding a
well founded structure satisfying $T$,'' and so it is
absolute between $V$ and $L$ by the Shoenfield absoluteness
theorem, provided the theory itself is in $L$.
Consequently, one can have transitive models of extremely
strong theories, including huge cardinals and more, without
ever leaving $L$.

Using this phenomenon, one may attack Steel's point 2 by
proposing the translation $\varphi\mapsto$``$\varphi$ holds
in a transitive model of set theory''. The believer in
large cardinals is usually happy to consider those
cardinals and other set-theoretic properties inside a
transitive model, and surely understanding how a particular
set-theoretic concept behaves inside a transitive model of
set theory is nearly the same as understanding how it
behaves in $V$. But the point is that with such a
relativization to transitive models, we get the same answer
in $L$ as in $V$, and so the assumption of $V=L$ did not
prevent us from asking about $\varphi$ in this way.

The large cardinal skeptics will object to Steel's point 1
on the grounds that the translation
$\varphi\mapsto\varphi^L$, if undertaken in the large
cardinal context in which Steel appears to intend it,
carries with it all the large cardinal baggage from such a
large cardinal $V$ to its corresponding $L$. After all, the
$L$ of a model with supercompact cardinals, for example,
will necessarily satisfy $\Con(\ZFC+\text{Woodin
cardinals})$ and much more, and so the strategy of working
in $V$ with large cardinals does not really offer a
translation for the believers in $V=L$ who are skeptical
about those large cardinals. Similarly, the relevant
comparison in point 3 would seem instead to be $\ZFC+V=L$
versus $\ZFC+\text{large cardinals}$, and here the simple
equivalence breaks down, precisely because the large
cardinals have consequences in $L$, such as the existence
of the corresponding shadow large cardinals in $L$, as well
as consistency statements that will be revealed in the
arithmetic truths of $L$. If \ZFC\ is consistent, then it
is of course incorrect to say that $\ZFC+V=L\proves\varphi$
if and only if $\ZFC+\text{large cardinals}\proves
\varphi^L$.

Pushing a bit harder on the examples above reveals that
every real of $V$ exists in a model $M$ of any desired
fixed finite fragment $\ZFC^*$ of $\ZFC$, plus $V=L$, which
furthermore is well-founded to any desired countable
ordinal height. That is, for every real $x$ and countable
ordinal $\alpha$, there is a model $M\satisfies\ZFC^*+V=L$,
whose well-founded initial segment exceeds $\alpha$ and
such that $x$ is among the reals of $M$. The proof of this
is that the statement is true inside $L$, and it has
complexity $\Pi^1_2$, adjoining a code for $\alpha$ as a
parameter if necessary, and hence true by the Shoenfield
absoluteness theorem. One can achieve full $\ZFC$ in the
model $M$ if $L$ has \ZFC\ models $L_\alpha$ unbounded in
$\omega_1^L$, which would hold if it had any uncountable
transitive model of \ZFC. Results in
\cite{HamkinsLinetskyReitz:PointwiseDefinableModelsOfSetTheory}
show moreover that one can even arrange that the model $M$
is pointwise definable, meaning that every set in $M$ is
definable there without parameters.

A striking example of the previous observation arises if
one supposes that $0^\sharp$ exists. Considering it as a
real, the previous remarks show that $0^\sharp$ exists
inside a pointwise definable model of $\ZFC+V=L$,
well-founded far beyond $\omega_1^L$. Thus, true $0^\sharp$
sits unrecognized, yet definable, inside a model of $V=L$
which is well-founded a long way. For a second striking
example, force to collapse $\omega_1$ to $\omega$ by a
generic filter $g$, which can be coded as a real. In
$V[g]$, there is $M\satisfies\ZFC+V=L$ with $g\in M$ and
$M$ well-founded beyond $\omega_1^V$. So $M$ thinks that
the generic real $g$ is constructible, constructed at some
(necessarily nonstandard) stage beyond $\omega_1^V$.

Another sense in which the axiom of constructibility is not
inherently limiting is provided by the observation that
every transitive model of set theory can in principle be
continued upward to a model of $V=L$. For a countable model
$M$, this follows from the previous observations about
reals, since we simply code $M$ with a real and then place
that real, and hence also $M$, into an $\omega$-model of
$\ZFC+V=L$, which in consequence will be well-founded
beyond the height of $M$. More generally, for any
transitive model, we may first collapse it to be countable
by forcing, and then carry out the previous argument in the
forcing extension. In this way, any transitive set $M$ can
in principle exist as a transitive set inside a model of
$\ZFC+V=L$. For all we know, our entire current universe,
large cardinals and all, is a countable transitive model
inside a much larger model of $V=L$.

In consequence, the axiom of constructibility $V=L$ begins
to shed its reputation as a limiting principle; perhaps it
will gain the reputation instead as a {\it rewarder of the
patient}, in the sense that if one waits long enough, then
one's sets become constructible. Much of the view of $V=L$
as inherently limiting stems from an implicit commitment to
the idea of an absolute background concept of ordinal and
the assumption that our current universe $V$ has all the
ordinals. That is, if we regard the current $L$ as already
containing all the ordinals, then the universe $V$ can only
grow sideways out of $L$, and with further sideways growth,
one never recovers $V=L$. But on the multiverse view, one
loses confidence in the idea of an absolute background
concept of the class of all ordinals, and we saw above how
easily the $V=L$ hypothesis can be recovered from any model
of set theory by growing upward as well as outward, rather
than only outward. As more ordinals arrive, we begin to see
new constructible sets that were formerly thought
non-constructible, and the constructible universe $L$ can
grow to encompass any given universe. In light of these
considerations, the models of $V=L$ do not seem to be
particularly impoverished.

\section{Multiverse Axioms}\label{Section.multiverseAxioms}

Let me turn now to a more speculative enterprise by
proposing a number of principles that express a fuller
vision of the multiverse perspective. These principles,
some of them provocative, should express the universe
existence properties we might expect to find in the
multiverse.

Before describing these principles, however, let me briefly
remark on the issue of formalism. Of course we do not
expect to be able to express universe existence principles
in the first-order language of set theory or even in the
second-order language of \Godel-Bernays or Kelly-Morse set
theory, since the entire point of the multiverse
perspective is that there may be other universes outside a
given one. In any of the usual formalizations of set
theory, we are of course at a loss to refer to these other
universes or to the objects in them. Thus, our foundational
ideas in the mathematics and philosophy of set theory
outstrip our formalism, mathematical issues become
philosophical, and set theory increasingly finds itself in
need of philosophical assistance. The difficulty is that
this philosophical need is arising in connection with some
of the most highly technical parts of the subject, such as
forcing, large cardinals and canonical inner models. So the
subject requires both technically minded philosophers and
philosophically minded mathematicians.

It turns out that many interesting questions about the
multiverse, such as those arising in the modal logic of
forcing or in set-theoretic geology (two topics explained
in the appendix of this article), can be formalized
entirely in the usual first-order language of set theory.
And for the issues in this realm, the work can be
understood as purely mathematical. I find it a fascinating
situation, where what might otherwise have been a purely
philosophical enterprise, the consideration and analysis of
alternative mathematical universes, becomes a fully
rigorous mathematical one. But for many other natural
questions about the multiverse, we seem unable adequately
to formalize them, despite our interest in the true
second-order or higher order issues. One fall-back method
is the {\it toy model perspective}, discussed at length in
\cite{FuchsHamkinsReitz:Set-theoreticGeology}, where one
considers a multiverse merely as a set of models in the
context of a larger \ZFC\ background. This approach is
adopted in
\cite{GitmanHamkins2010:NaturalModelOfMultiverseAxioms} and
also in \cite{Friedman2006:InternalConsistencyAndIMH,
FriedmanWelchWoodin2008:ConsistencyOfIMH} in the context of
the inner model hypothesis. It is precisely the toy model
perspective that underlies the traditional countable
transitive model approach to forcing. The toy model
perspective can ultimately be unsatisfying, however, since
it is of course in each case not the toy model in which we
are interested, but rather the fully grown-up universe.
Allow me, then, to state the multiverse axioms as
unformalized universe existence assertions about what we
expect of the genuine full multiverse.

The background idea of the multiverse, of course, is that
there should be a large collection of universes, each a
model of (some kind of) set theory. There seems to be no
reason to restrict inclusion only to \ZFC\ models, as we
can include models of weaker theories $\ZF$, $\ZF^-$, $\KP$
and so on, perhaps even down to second order number theory,
as this is set-theoretic in a sense. At the same time,
there is no reason to consider all universes in the
multiverse equally, and we may simply be more interested in
the parts of the multiverse consisting of universes
satisfying very strong theories, such as \ZFC\ plus large
cardinals. The point is that there is little need to draw
sharp boundaries as to what counts as a set-theoretic
universe, and we may easily regard some universes as more
set-theoretic than others.

Nevertheless, we have an expansive vision for the
multiverse. Explaining set-theorist's resistance to $V=L$
as an axiom, Maddy \cite{Maddy1998:V=LAndMaximize} has
pointed to the MAXIMIZE maxim---asserting that there should
be no undue limitations on the sets that might exist---for
it would inflict an artificial limitation to insist that
every set is constructible from ordinals. Here, we may
follow a multiverse analogue of MAXIMIZE, by placing no
undue limitations on what universes might exist in the
multiverse. This is simply a higher-order analogue of the
same motivation underlying set-theorists' ever more
expansive vision of set theory. We want to consider that
the multiverse is as big as we can imagine. At any time, we
are living inside one of the universes, referred to as $V$
and thought of as the current universe, but by various
means, sometimes metamathematical, we may be able to move
around in the multiverse. We do not expect individual
worlds to access the whole multiverse. Thus, we do not want
to allow set construction principles in any one universe
that involve quantifying over all universes, unless this
quantification can be shown to reduce to first order
quantification over sets in that universe.

Let me begin with the most basic universe existence
principles. It goes without saying that the view holds
first, that there is a universe. Going beyond this, the
following principle asserts that inner models of a given
universe, and more generally interpreted models, are also
universes in their own right.

\begin{proclamation}[Realizability Principle.]
For any universe $V$, if $W$ is a model of set theory and
definable or interpreted in $V$, then $W$ is a universe.
\end{proclamation}

The principle expresses the multiverse capacity to regard
constructed universes as real. In algebra, when one group
is constructed from another or when a field is constructed
as a quotient of a ring, these new structures are taken to
be perfectly good groups and fields. There is no objection
in such a case announcing ``{\it but that isn't the {\it
real} $+$},'' analogous to what one might hear about a
nonstandard or disfavored set-membership relation: ``{\it
that isn't the {\it real} $\in$}.'' Such remarks in set
theory are especially curious given that the nature and
even the existence of the intended absolute background
notion is so murky. On the multiverse view, if someone
constructs a model of set theory, then we simply wonder,
``what would it be like to live in that universe?'' Jumping
inside, we want to regard it as the real universe. The
realizability principle allows this, and encompasses all
the usual definable inner models, such as $L$, $\HOD$,
$L(\R)$, as well as ultrapowers and many others, including
models that are sets in $V$, such as $J_\beta$, $V_\alpha$
and $H_\kappa$, and so on. In at least a weak sense, the
realizability principle should be unproblematic even for
those with the universe view, since even if you have an
absolute background universe, then the inner models have a
real existence in at least some sense.

The next principle, in contrast, posits the existence of
certain outer models and therefore engages with the main
disagreement between the universe and multiverse views. I
regard the following principle as a minimal foundational
commitment for the multiverse view.

\begin{proclamation}[Forcing Extension Principle.]
 For any
universe $V$ and any forcing notion $\P$ in $V$, there is a
forcing extension $V[G]$, where $G\of\P$ is $V$-generic.
\end{proclamation}

\noindent Naturally, we also expect many non-forcing
extensions as well.

Next I consider a multiverse principle I refer to as the
{\it reflection axiom}. Every model of set theory imagines
itself as very tall, containing all the ordinals. But of
course, a model contains only all of its {\it own}
ordinals, and we know very well that it can be mistaken
about how high these reach, since other models may have
many more ordinals. We have a rich experience of this
phenomenon with models of set theory, and many arguments
rely precisely on this point, of a model thinking it has
all the ordinals, when actually it is very small, even
countable. Thus, on the multiverse perspective, we should
expect that every universe can be extended to a much taller
universe. In addition, many philosophers have emphasized
the importance of reflection principles in set theory,
using them to justify increasingly strong large cardinal
axioms and to fulfill the ideas at the dawn of set theory
of the unknowability of the full universe of sets. Thus, I
propose the following multiverse principle:

\begin{proclamation}[Reflection Axiom.]
For every universe $V$, there is a much taller universe $W$
with an ordinal $\theta$ for which $V\precsim W_\theta\prec
W$.
\end{proclamation}

The toy model version of this axioms is easily seen to be a
consequence of the compactness theorem, applied to the
elementary diagram of $V$ and its relativization to
$W_\theta$, using a new ordinal constant $\theta$. But
here, in the full-blown multiverse context, the principle
asserts that no universe is correct about the height of the
ordinals, and every universe looks like an initial segment
of a much taller universe having the same truths.

Next, consider a somewhat more provocative principle. Every
model of set theory imagines itself as enormous, containing
the reals $\R$, all the ordinals, and so on, and it
certainly believes itself to be uncountable. But of course,
we know that models can be seriously mistaken about this.
The Lowenheim-Skolem theorem, of course, explains how there
can be countable models of set theory, and these models do
not realize there is a much larger universe surrounding
them. Countability after all is a notion that depends on
the set-theoretic background, since the very same set can
be seen either as countable or uncountable in different
set-theoretic contexts. Indeed, an elementary forcing
argument shows that any set in any model of set theory can
be made countable in a suitable forcing extension of that
universe. On the multiverse perspective, the notion of an
absolute standard of countability falls away. I therefore
propose the following principle:

\begin{proclamation}[Countability Principle.] Every universe
$V$ is countable from the perspective of another, better
universe $W$.
\end{proclamation}

If one merely asserts that $V$ is an element of another
universe $W$, then by the Forcing Extension principle it
will be countable in the forcing extension $W[G]$ that
collapses suitable cardinals to $\omega$, and so there is a
measure of redundancy in these multiverse axioms, here and
in several other instances.

Let me consider next what may be the most provocative
principle. Set-theorists have traditionally emphasized the
well-founded models of set theory, and the Mostowski
collapses of these models are precisely the transitive
models of set theory, universally agreed to be highly
desirable models. The concept of well-foundedness, however,
depends on the set-theoretic background, for different
models of set theory can disagree on whether a structure is
well-founded. Every model of set theory finds its own
ordinals to be well-founded, but we know that models can be
incorrect about this, and one model $W$ can be realized as
ill-founded by another, better model. This issue arises, as
I emphasized earlier, even with the natural numbers. Every
model of set theory finds its own natural numbers $\N$ to
be `the standard' model of arithmetic, but different models
of set theory can have different (non-isomorphic)
`standard' models of arithmetic. The basic situation is
that the $\N$ of one model of set theory can often be
viewed as non-standard or even uncountable from another
model. Indeed, every set theoretic argument can take place
in a model, which from the inside appears to be totally
fine, but actually, the model is seen to be ill-founded
from the external perspective of another, better model.
Under the universe view, this problem terminates in the
absolute set-theoretic background universe, which provides
an accompanying absolute standard of well-foundedness. But
the multiverse view allows for many different set-theoretic
backgrounds, with varying concepts of the well-founded, and
there seems to be no reason to support an absolute notion
of well-foundedness. Thus, we seem to be led to the
following principle:

\begin{proclamation}[Well-foundedness Mirage.] Every
universe $V$ is ill-founded from the perspective of
another, better universe.
\end{proclamation}

This principle appears to be abhorrent to many set
theorists---and so I hope to have succeeded in my goal to
be provocative---but I have a serious purpose in mind with
it. The point is that we know very well that the concept of
well-foundedness depends on the set-theoretic background,
and so how can we accept philosophical talk focussing on
well-founded models, using this term as though it had an
absolute meaning independent of any set-theoretic
background? And once we give up on an absolute
set-theoretic background, then any grounding we may have
had for an absolute background concept of well-foundedness
falls away. On the multiverse view, therefore, for any
universe there may be other universes with an improved
concept of well-foundedness, and others with still better
concepts, and we may apply this idea even to the natural
numbers of the models.

Consider next the reverse embedding axiom. When we are
living in a universe $V$ and have an embedding $j:V\to M$,
such as an ultrapower embedding, then we know very well how
to iterate this map forwards
$$V\to M\to M_2\to M_3\to\cdots$$
The later models $M_n$ have their iterated version of the
embedding $j_n:M_n\to M_{n+1}$, and do not realize that
these maps have already been iterated many times. On the
multiverse view, every embedding is like this, having
already been iterated an enormous number of times.

\begin{proclamation}[Reverse Embedding Axiom.]
For every universe $V$
    and every embedding $j:V\to M$ in
    $V$, there is a universe $W$ and embedding $h$
  $$W\rTo^h V\rTo^j M$$
   such that $j$ is the iterate of $h$.
\end{proclamation}

We saw earlier how every countable transitive model of set
theory can be continued to a model of $V=L$. On the
multiverse view, we expect every universe to be similarly
absorbed into $L$, so that every universe is a countable
transitive model inside a much larger $L$. In this sense,
$L$ is a rewarder of the patient, for given sufficient
ordinals your favorite set and indeed, your whole universe,
will become both countable and constructible.

\begin{proclamation}[Absorption into $L$.]
Every universe $V$ is a countable transitive model in
another universe $W$ satisfying $V=L$.
\end{proclamation}

I have described a number of multiverse axioms, sketching
out a vision of the universe existence principles we expect
to find in the multiverse. A natural question, of course,
is whether this vision is coherent. One way to mathematize
this coherence question is to consider the multiverse
axioms in their toy model interpretation: can one have a
set collection of models of \ZFC\ that satisfy all of the
closure properties expressed by the multiverse axioms
above? And for this toy model version of the question, the
answer is yes. Victoria Gitman and I proved in
\cite{GitmanHamkins2010:NaturalModelOfMultiverseAxioms}
that the collection of countable computably-saturated
models of \ZFC\ satisfies all of the multiverse axioms
above. With acknowledgement for the limitations of the toy
model perspective as a guide to the full, true,
higher-order multiverse, we take this as evidence that the
multiverse vision expressed by the multiverse axioms above
is coherent.

In conclusion, let me remark that this multiverse vision,
in contrast to the universe view with which we began this
article, fosters an attitude that what set theory is about
is the exploration of the extensive range of set-theoretic
possibilities. And this is an attitude that will lead any
young researcher to many interesting set-theoretic
questions and, I hope, fruitful answers.

\section{Appendix: multiverse-Inspired
Mathematics}\label{Appendix}

The mathematician's measure of a philosophical position may
be the value of the mathematics to which it leads. Thus,
the philosophical debate here may be a proxy for: where
should set theory go? Which mathematical questions should
it consider? The multiverse view in set theory leads one to
consider how the various models of set theory interact, or
how a particular world sits in the multiverse. This would
include such mainstream topics as absoluteness,
independence, forcing axioms, indestructibility and even
large cardinal axioms, which concern the relation of the
universe $V$ to certain inner models $M$, such as an
ultrapower. A few recent research efforts, however, have
exhibited the multiverse perspective more explicitly, and I
would like briefly to describe two such projects in which I
have been involved, namely, the modal logic of forcing and
set-theoretic geology. This appendix is adapted from
\cite{HamkinsLoewe2008:TheModalLogicOfForcing,
FuchsHamkinsReitz:Set-theoreticGeology} and the surveys
\cite{Hamkins2009:SomeSecondOrderSetTheory,
Hamkins:Themultiverse:ANaturalContext}.

The modal logic of forcing (see
\cite{HamkinsLoewe2008:TheModalLogicOfForcing} for a full
account) is concerned with the most general principles in
set theory relating forcing and truth. Although we did not
begin with the intention to introduce modal logic into set
theory, these principles are most naturally expressed in
modal logic, for the collection of forcing extensions forms
a robust Kripke model of possible worlds. Each such
universe serves as an entire mathematical world and every
ground model has a degree of access, via names and the
forcing relation, to the objects and truths of its forcing
extensions. Thus, this topic engages explicitly with a
multiverse perspective. The natural forcing modalities
(introduced in \cite{Hamkins2003:MaximalityPrinciple})
assert that a sentence $\varphi$ is {\df possible} or {\df
forceable}, written {\df $\possible\varphi$}, when it holds
in a forcing extension and {\df necessary}, written {\df
$\necessary\varphi$}, when it holds in all forcing
extensions. Thus, $\varphi$ is forceably necessary, if it
is forceable that $\varphi$ is necessary; in other words,
if there is a forcing extension in which $\varphi$ holds
and continues to hold in all further extensions.

The initial inquiry centered on the question: Could the
universe be completed with respect to what is forceably
necessary? That is, could there be a model over which any
possibly necessary assertion is already necessary? The {\it
Maximality Principle} is the scheme expressing exactly
that:
$$\possible\necessary\varphi\implies\necessary\varphi$$
This principle asserts that the universe has been completed
with respect to forcing in the sense that everything
permanent achievable by forcing has already been achieved.
It turns out to be equiconsistent with \ZFC.

\begin{theorem}[{Hamkins \cite{Hamkins2003:MaximalityPrinciple}, also Stavi \&
V\"a\"an\"anen \cite{StaviVaananen:ReflectionPrinciples},
independently}] The Maximality Principle is relatively
consistent with \ZFC.
\end{theorem}

Another way to state the theorem is that if \ZFC\ is
consistent, then there is a model of \ZFC\ in which the
implication
$\possible\necessary\varphi\implies\necessary\varphi$ is a
valid principle of forcing, in the sense that it holds
under the forcing interpretation for all sentences
$\varphi$ in the language of set theory. This led to a deep
question:

\begin{question}
What are the generally correct principles of forcing? Which
modal principles are valid in all models of set theory?
Which can be valid?
\label{Question.ValidPrinciplesOfForcing?}
\end{question}

This question and its solution are deeply engaged with the
multiverse perspective. L\"owe and I were inspired by
Solovay's \cite{Solovay1976:ProvabilityModalLogic} modal
analysis of provability, in which he discovered the
corresponding valid principles of provability, and we aimed
to do for forceability what Solovay did for provability.

To get started, the following principles are trivially
valid principles of forcing:

$$\begin{array}{rl}
 \axiomf{K} & \necessary(\varphi\implies\psi)\implies(\necessary\varphi\implies\necessary\psi)\\
 \axiomf{Dual} & \necessary\neg\varphi\iff\neg\possible\varphi\\
 \axiomf{S} & \necessary\varphi\implies\varphi\\
 \axiomf{4} & \necessary\varphi\implies\necessary\necessary\varphi\\
 \axiomf{.2} & \possible\necessary\varphi\implies\necessary\possible\varphi\\
\end{array}$$

For example, \axiomf{K} asserts that if an implication
holds in every forcing extension and also its hypothesis,
then so does its conclusion. Axioms \axiomf{Dual},
\axiomf{S} and \axiomf{4} are similarly easy. Axiom .2
asserts that every possibly necessary statement is
necessarily possible, and this is true because if $\P$
forces $\varphi$ to be necessary, and $\Q$ is some other
forcing notion, then forcing with $\P$ over $V^\Q$ is
equivalent by product forcing to forcing with $\Q$ over
$V^\P$, and so $\varphi$ is possible over $V^\Q$, as
desired. Since the assertions above axiomatize the modal
theory known as \axiomf{S4.2} (and because the \ZFC\
provable validities are closed under necessitation), it
follows that any \axiomf{S4.2} modal assertion is a valid
principle of forcing.

Going beyond \axiomf{S4.2}, the following is a sample of
modal assertions not derivable in \axiomf{S4.2}---any modal
logic text will contain hundreds more---each of which is
important in some variant modal context. For example, the
\axiomf{L\"ob} axiom figures in Solovay's solution to
provability logic, and expresses the provability content of
L\"ob's theorem. Does the reader think that any of these
express valid principles of forcing?
$$\begin{array}{rl}
 \axiomf{5} & \possible\necessary\varphi\implies\varphi \\
 \axiomf{M} & \necessary\possible\varphi\implies\possible\necessary\varphi \\
 \axiomf{W5} & \possible\necessary\varphi\implies(\varphi\implies\necessary\varphi) \\
 \axiomf{.3} & \possible\varphi\wedge\possible\psi\implies(\possible(\varphi\wedge\possible\psi)
             \vee\possible(\varphi\wedge\psi)\vee\possible(\psi\wedge\possible\varphi))\\

 \axiomf{Dm} & \necessary(\necessary(\varphi\implies\necessary\varphi)\implies\varphi)
               \implies(\possible\necessary\varphi\implies\varphi)\\
 \axiomf{Grz} & \necessary(\necessary(\varphi\implies\necessary\varphi)\implies\varphi)\implies\varphi \\
 \axiomf{\hbox{L\"ob}} & \necessary(\necessary\varphi\implies\varphi)\implies\necessary\varphi \\
 \axiomf{H} & \varphi\implies\necessary(\possible\varphi\implies\varphi)\\
\end{array}$$
It is a fun forcing exercise to show that actually, none of
these assertions is valid for the forcing interpretation in
all models of \ZFC, although some are valid in some but not
all models and others are never valid. The resulting modal
theories can be arranged by strength into the following
diagram, which originally focused our attention on
\theoryf{S4.2} as a possible answer to question
\ref{Question.ValidPrinciplesOfForcing?}, since this theory
is valid but all stronger theories in the diagram fail in
some models of set theory.
$$\quad\hbox{\tiny\vbox{\hbox{$\begin{array}{rcl}
 \theoryf{S5}    &=& \theoryf{S4}+5  \\
 \theoryf{S4W5}&=& \theoryf{S4}+\axiomf{W5}\\
 \theoryf{S4.3}&=& \theoryf{S4}+\axiomf{.3}\\
 \theoryf{S4.2.1}&=& \theoryf{S4}+\axiomf{.2}+\axiomf{M}\\
 \theoryf{S4.2}&=& \theoryf{S4}+\axiomf{.2}\\
 \theoryf{S4.1}&=& \theoryf{S4}+\axiomf{M}\\
 \theoryf{S4}     &=& \theoryf{K4}+\axiomf{S}\\
 \theoryf{Dm.2}&=& \S{4.2}+\axiomf{Dm}\\
 \theoryf{Dm}&=& \S4+\axiomf{Dm}\\
 \theoryf{Grz}&=& \theoryf{K}+\axiomf{Grz}\\
 \theoryf{GL}&=& \theoryf{K4}+\axiomf{\hbox{L\"ob}}\\
 \theoryf{K4H}&=& \theoryf{K4}+\axiomf{H}\\   % note that this is the *new* H axiom
 \theoryf{K4} &=& \theoryf{K}+4\\
 \theoryf{K} &=& \axiomf{K}+\axiomf{Dual}\\
\end{array}$}%
}}%
\hbox{\tiny
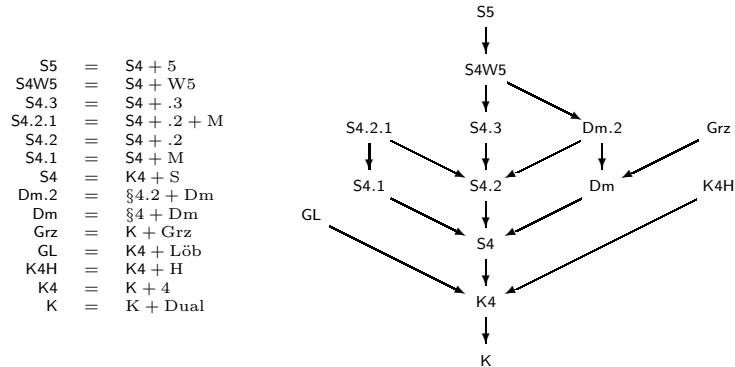
\begin{diagram}[height=1.5em]
            &           &           &           &  {\theoryf{S5}}      &           &             &            &        \\
            &           &           &           &  \dTo     &           &             &            &        \\
            &           &           &           &  {\theoryf{S4W5}}  &           &             &            &        \\
            &           &           &           &  \dTo     &  \rdTo    &             &            &        \\
            &           & {\theoryf{S4.2.1}} &           & {\theoryf{S4.3}}   &           & {\mathsf{Dm.2}}&            & {\mathsf{Grz}} \\
            &           & \dTo      & \rdTo     &  \dTo     &  \ldTo    &   \dTo      & \ldTo      &           \\
            &           & {\theoryf{S4.1}}   &           & {\theoryf{S4.2}}   &           & {\mathsf{Dm}}  &            & {\mathsf{K4H}}  \\
            &{\mathsf{GL}}&           & \rdTo     &  \dTo     &  \ldTo    &             & \ldTo(4,4) &  \\
            &           & \rdTo(3,3)&           &{\theoryf{S4}}      &           &             &            &        \\
            &           &           &           & \dTo      &           &             &            &       \\
            &           &           &           &{\mathsf{K4}}&           &             &            &        \\
            &           &           &           &  \dTo     &           &             &            &        \\
            &           &           &           &{\mathsf{K}} &           &             &            &        \\
\end{diagram}}$$
L\"owe and I were very pleased to prove that this answer
was indeed correct:

\begin{theorem}[Hamkins,L\"owe
\cite{HamkinsLoewe2008:TheModalLogicOfForcing}] If \ZFC\ is
consistent, then the \ZFC-provably valid principles of
forcing are exactly
\theoryf{S4.2}.\label{Theorem.ModalLogicOfForcingIsS4.2}
\end{theorem}

We observed above that \theoryf{S4.2} is valid; the
difficult part of the theorem was to show that nothing else
is valid. If $\theoryf{S4.2}\not\proves\varphi$, then we
must provide a model of set theory and set-theoretic
assertions $\psi_i$ such that
$\varphi(\psi_0,\ldots,\psi_n)$ fails in that model. But
what set-theoretic content should these assertions express?
It seems completely mysterious at first. I would like to
mention a few of the ideas from the proof. Two attractively
simple concepts turn out to be key. Specifically, a
statement $\varphi$ of set theory is said to be a {\df
switch} if both $\varphi$ and $\neg\varphi$ are necessarily
possible. That is, a switch is a statement whose truth
value can always be turned on or off by further forcing. In
contrast, a statement is a {\df button} if it is
(necessarily) possibly necessary. These are the statements
that can be forced true in such a way that they remain true
in all further forcing extensions. The idea here is that
once you push a button, you cannot unpush it. The
Maximality Principle, for example, asserts that every
button has already been pushed. Although buttons and
switches may appear at first to be very special kinds of
statements, nevertheless every statement in set theory is
either a button, a switch, or the negation of a button.
(After all, if you can't always switch $\varphi$ on and
off, then it will either get stuck on or stuck off, and
product forcing shows these possibilities to be mutually
exclusive.) A family of buttons and switches is {\df
independent}, if the buttons are not yet pushed and
(necessarily) each of the buttons and switches can be
controlled without affecting the others. Under $V=L$, it
turns out that there is an infinite independent family of
buttons and switches.

The proof of Theorem
\ref{Theorem.ModalLogicOfForcingIsS4.2} rests in part on a
detailed understanding of the modal logic \theoryf{S4.2}
and its complete sets of Kripke frames. A (propositional)
{\df Kripke model} is a collection of propositional {\df
worlds} (essentially a truth table row, assigning
propositional variables to true and false), with an
underlying accessibility relation called the {\df frame}. A
statement is possible or necessary at a world, accordingly
as it is true in some or all accessible worlds,
respectively. Every Kripke model built on a frame that is a
directed partial pre-order will satisfy the \theoryf{S4.2}
axioms of modal logic, and in fact the finite directed
partial pre-orders are {\df complete} for \theoryf{S4.2} in
the sense that the statements true in all Kripke models
built on such frames are exactly the statements provable
from \theoryf{S4.2}. An improved version of this, proved in
\cite{HamkinsLoewe2008:TheModalLogicOfForcing}, is that the
finite pre-lattices, and even the finite pre-Boolean
algebras, are complete for \theoryf{S4.2}. The following
lemma, a central technical claim at the core of
\cite{HamkinsLoewe2008:TheModalLogicOfForcing}, shows that
any model of set theory with an independent family of
buttons and switches is able to simulate any given Kripke
model built on a finite pre-lattice frame.

\begin{lemma}\label{Lemma.Simulation} If $W\satisfies\ZFC$
has sufficient independent buttons and switches, then for
any Kripke model $M$ on a finite pre-lattice frame, any
$w\in M$, there is a translation of the propositional
variables $p_i\mapsto\psi_i$ to set-theoretic assertions
$\psi_i$, such that for any modal assertion
$\varphi(p_1,\ldots,p_n)$:
$$(M,w)\satisfies\varphi(p_1,\ldots,p_n)\quad\iff\quad W\satisfies\varphi(\psi_1,\ldots,\psi_n).$$\smallskip
Each $\psi_i$ is a Boolean combination of the buttons and
switches.
\end{lemma}

Consequently, if $\theoryf{S4.2}\not\proves\varphi$, then
since we proved that there is a Kripke model $M$ built on a
finite pre-lattice frame in which $\varphi$ fails, it
follows that in any model of set theory $W$ having
independent buttons and switches, which we proved exist,
the corresponding assertion $\varphi(\psi_1,\ldots,\psi_n)$
fails. This exactly shows that $\varphi$ is not a provably
valid principle of forcing, as desired to prove Theorem
\ref{Theorem.ModalLogicOfForcingIsS4.2}. The proof is
effective in the sense that if
$\theoryf{S4.2}\not\proves\varphi$, then we are able
explicitly to provide a model $W\satisfies\ZFC$ and the
particular set-theoretic substitution instance
$\varphi(\psi_1,\ldots,\psi_n)$ which fails in $W$. To
answer our earlier query about what the content of the
sentences $\psi_i$ might be, the answer is that they are
certain complex Boolean combinations of the independent
buttons and switches relating to the Jankov-Fine formula
describing the underlying frame of a propositional Kripke
model in which $\varphi$ fails. The buttons and switches
can be combined so as to fail under the forcing
interpretation in just exactly the way that the modal
assertion $\varphi$ fails in a finite Kripke model.

Although Theorem \ref{Theorem.ModalLogicOfForcingIsS4.2}
tells us what are the \ZFC-provably valid principles of
forcing, it does not tell us that all models of \ZFC\
exhibit only those validities. Indeed, we know that this
isn't the case, because we know there are models of the
Maximality Principle, for which the modal theory
\theoryf{S5} is valid, and this is strictly stronger than
\theoryf{S4.2}. So different models of set theory can
exhibit different valid principles of forcing. The proof of
Theorem \ref{Theorem.ModalLogicOfForcingIsS4.2} adapts to
show, however, that the valid principles of forcing over
any fixed model of set theory lies between $\theoryf{S4.2}$
and $\theoryf{S5}$. This result is optimal because both of
these endpoints occur: the theory $\theoryf{S4.2}$ is
realized in any model of $V=L$, and $\theoryf{S5}$ is
realized in any model of the Maximality Principle. This
result also shows that the Maximality Principle is the
strongest possible axiom of this type, expressible in this
modal language as valid principles of forcing. A number of
questions remain open, such as: Is there a model of \ZFC\
whose valid principles of forcing form a theory other than
\theoryf{S4.2} or \theoryf{S5}? If $\varphi$ is valid in
$W$, is it valid in all extensions of $W$? Can a model of
\ZFC\ have an unpushed button, but not two independent
buttons?

Let me turn now to the topic of set-theoretic geology,
another topic engaged with the multiverse perspective (see
\cite{FuchsHamkinsReitz:Set-theoreticGeology} for a full
introductory account). Although forcing is ordinarily
viewed as a method of constructing {\it outer} as opposed
to {\it inner} models of set theory---one begins with a
ground model and builds a forcing extension of it---a
simple switch in perspective allows us to use forcing to
describe inner models as well. The idea is simply to
consider things from the perspective of the forcing
extension rather than the ground model and to look downward
from the universe $V$ for how it may have arisen by
forcing. Given the set-theoretic universe $V$, we search
for the possible ground models $W\of V$ such that there is
a $W$-generic filter $G\of\P\in W$ such that $V=W[G]$. Such
a perspective quickly leads one to look for deeper and
deeper grounds, burrowing down to what we call bedrock
models and deeper still, to what we call the mantle and the
outer core, for the topic of set-theoretic geology.

Haim Gaifman has pointed out (with humor) that the geologic
terminology takes an implicit position on the question of
whether mathematics is a natural process or a human
activity, for otherwise we might call it set-theoretic {\it
archeology}, aiming to burrow down and uncover the remains
of ancient set-theoretic civilizations.

Set-theoretic geology rests fundamentally on the following
theorem, a surprisingly recent result considering the
fundamental nature of the question it answers. Laver's
proof of this theorem make use of my work on the
approximation and covering properties
\cite{Hamkins2003:ExtensionsWithApproximationAndCoverProperties}.

\begin{theorem}[{Laver
\cite{Laver2007:CertainVeryLargeCardinalsNotCreated},
independently Woodin
\cite{Woodin2004:RecentDevelopmentsOnCH}}] Every model of
set theory $V\satisfies\ZFC$ is a definable class in all of
its set forcing extensions $V[G]$, using parameters in
$V$.\label{Theorem.GroundIsDefinable}
\end{theorem}

After learning of this theorem, Jonas Reitz and I
introduced the following hypothesis:

\begin{definition}[{Hamkins, Reitz
\cite{Hamkins2005:TheGroundAxiom,Reitz2006:Dissertation,Reitz2007:TheGroundAxiom}}]\rm
The ground axiom \GA\ is the assertion that the universe is
not obtained by nontrivial set forcing over any inner
model.
\end{definition}

Because of the quantification over ground models, this
assertion appears at first to be second order, but in fact
the ground axiom is expressible by a first order statement
in the language of set theory (see
\cite{Reitz2006:Dissertation,Reitz2007:TheGroundAxiom}).
And although the ground axiom holds in many canonical
models of set theory, such as $L$, $L[0^\sharp]$, $L[\mu]$
and many instances of $K$, none of the regularity features
of these models are consequences of \GA, for Reitz proved
that every \ZFC\ model has an extension, preserving any
desired $V_\theta$ (and mild in the sense that every new
set is generic for set forcing), which is a model of \ZFC\
plus the ground axiom. Reitz's method obtained \GA\ by
forcing very strong versions of $V=\HOD$, but in a
three-generation collaboration between Reitz, Woodin and
myself
\cite{HamkinsReitzWoodin2008:TheGroundAxiomAndV=HOD}, we
proved that every model of set theory has an extension
which is a model of \GA\ plus $V\neq\HOD$.

A {\df ground} of the universe $V$ is a transitive class
$W\satisfies\ZFC$ over which $V$ was obtained by forcing,
so that $V=W[G]$ for some $W$-generic $G\of\P\in W$. A
ground $W$ is a {\df bedrock} of $V$ if it is a minimal
ground, that is, if there is no deeper ground inside $W$.
Equivalently, $W$ is a bedrock of $V$ if it is a ground of
$V$ and satisfies the ground axiom.  It is an open question
whether a model can have more than one bedrock model.
However, Reitz \cite{Reitz2006:Dissertation} constructs a
{\df bottomless} model of \ZFC, one having no bedrock at
all. Such a model has many grounds, over which it is
obtained by forcing, but all of these grounds are
themselves obtained by forcing over still deeper grounds,
and so on without end. The principal new set-theoretic
concept in geology is the following:

\begin{definition}\rm The {\df mantle} $\Mantle$ is the
intersection of all grounds.
\end{definition}

The analysis at this point engages with an interesting
philosophical position, the view I call {\it ancient
paradise}. This position holds that there is a highly
regular core underlying the universe of set theory, an
inner model obscured over the eons by the accumulating
layers of debris heaped up by innumerable forcing
constructions since the beginning of time. If we could
sweep the accumulated material away, on this view, then we
should find an ancient paradise. The mantle, of course,
wipes away an entire strata of forcing, and so if one
subscribes to the ancient paradise view, then it seems one
should expect the mantle to exhibit some of these highly
regular features. Although I find the position highly
appealing, our initial main theorem appears to provide
strong evidence against it.

\begin{theorem}[{Fuchs, Hamkins, Reitz
\cite{FuchsHamkinsReitz:Set-theoreticGeology}}] Every model
of \ZFC\ is the mantle of another model of
\ZFC.\label{Theorem.EveryModelIsAMantle}
\end{theorem}

The theorem shows that by sweeping away the accumulated
sands of forcing, what we find is not a highly regular
ancient core, but rather: an arbitrary model of set theory.
In particular, since every model of set theory is the
mantle of another model, one can prove nothing special
about the mantle, and certainly, it need be no ancient
paradise.

The fact that the mantle $\Mantle$ is first-order definable
follows from the following general fact, which reduces
second-order quantification over grounds to first order
quantification over indices.

\begin{theorem}
There is a parameterized family of classes
$W_r$ such that
\begin{enumerate}
 \item Every $W_r$ is a ground of $V$ and $r\in W_r$.
 \item Every ground of $V$ is $W_r$ for some $r$.
 \item The relation ``$x\in W_r$'' is first order.
\end{enumerate}
\end{theorem}

For example, the ground axiom is first order expressible as
$\forall r\, W_r=V$; the model $W_r$ is a bedrock if and
only if $\forall s\, (W_s\of W_r\implies W_s=W_r)$; and the
mantle is defined by $\Mantle=\set{x\st \forall r\, (x\in
W_r)}$, and these are all first order assertions. When
looking downward at the various grounds, it is very natural
to inquire whether one can fruitfully intersect them. For
example, the question of whether there can be distinct
bedrock models in the universe is of course related to the
question of whether there is a ground in their
intersection. Let us define that the grounds are {\df
downward directed} if for every $r$ and $s$ there is $t$
such that $W_t\of W_r\intersect W_s$. In this case, we say
that the downward directed grounds hypothesis (\DDG) holds.
The strong \DDG\ asserts that the grounds are downward
set-directed, that is, that for every set $A$ there is $t$
with $W_t\of \Intersect_{r\in A}W_r$. These hypotheses
concern fundamental questions about the structure of ground
models of the universe---fundamental properties of the
multiverse---and these seem to be foundational questions
about the nature and power of forcing. In every model for
which we can determine the answer, the strong \DDG\ holds,
although it is not known to be universal. In any case, when
the grounds are directed, the mantle is well behaved.

\begin{theorem}\label{Theorem.IfDirectedMantleHasZF}\
\begin{enumerate}
 \item If the grounds are downward directed, then the
     mantle is constant across the grounds and a model
     of \ZF.
 \item If the grounds are downward set-directed, then
     the mantle is a model of \ZFC.
\end{enumerate}
\end{theorem}

It is natural to consider how the mantle is affected by
forcing, and since every ground of $V$ is a ground of any
forcing extension of $V$, it follows that the mantle of any
forcing extension is contained in the mantle of the ground
model. In the limit of this process, we arrive at:

\begin{definition}\rm The {\df generic mantle} of a model of set
theory $V$ is the intersection of all ground models of all
set forcing extensions of $V$.
\end{definition}

The generic mantle has proved in several ways to be a more
robust version of the mantle. For example, for any model of
\ZFC, the generic mantle $\gMantle$ is always a model at
least of \ZF, without need for any directedness hypothesis.
If the \DDG\ holds in all forcing extensions, then in fact
the mantle and the generic mantle are the same. If the
strong \DDG\ also holds, then the generic mantle is a model
of \ZFC.

Set-theoretic geology is naturally carried out in a context
that includes all the forcing extensions of the universe,
all the grounds of these extensions, all forcing extensions
of these resulting grounds and so on. The {\df generic
multiverse} of a model of set theory, introduced by Woodin
\cite{Woodin2004:RecentDevelopmentsOnCH}, is the smallest
family of models of set theory containing that model and
closed under both forcing extensions and grounds. There are
numerous philosophical motivations to study to the generic
multiverse. Indeed, Woodin introduced it specifically in
order to criticize a certain multiverse view of truth,
namely, truth as true in every model of the generic
multiverse. Although I do not hold such a view of truth,
nevertheless I want to investigate the fundamental features
of the generic multiverse, a task I place at the foundation
of any deep understanding of forcing. Surely the generic
multiverse is the most natural and illuminating background
context for the project of set-theoretic geology.

If the \DDG\ holds in all forcing extensions, then the
grounds are dense below the generic multiverse, and so not
only does $\Mantle=\gMantle$, but the generic multiverse is
exhausted by the ground extensions $W_r[G]$, the forcing
extensions of the grounds. In this case, any two models in
the generic multiverse have a common ground, and are
therefore at most two steps apart on a generic multiverse
path. There is a toy model example, however, showing that
the generic grounds---the ground models of the forcing
extensions---may not necessarily exhaust the generic
multiverse.

The generic mantle $\gMantle$ is a parameter-free uniformly
definable class model of \ZF, containing all ordinals and
invariant by forcing. Because it is invariant by forcing,
the generic mantle $\gMantle$ is constant across the
multiverse. In fact, it follows that the generic mantle
$\gMantle$ is precisely the intersection of the generic
multiverse. Therefore, it is the largest forcing-invariant
class. Surely this makes it a highly canonical class in set
theory; it should be the focus of intense study. Whenever
anyone constructs a new model of set theory, I should like
to know what is the mantle and generic mantle of this
model.

The {\df generic $\HOD$}, introduced by Fuchs in an attempt
to identify a very large canonical forcing-invariant class,
is the intersection of all $\HOD$s of all forcing
extensions $V[G]$.$$\gHOD=\Intersect_G \HOD^{V[G]}$$ This
class is constant across the generic multiverse. Since the
$\HOD$s of all the forcing extensions are downward
set-directed, it follows that $\gHOD$ is locally realized
and $\gHOD\satisfies\ZFC$. The following diagram shows the
basic inclusion relations.

\begin{diagram}[height=2em]
\HOD & \\
\cup \\   % improve this symbol, for vertical inclusion
\gHOD & \of & \gMantle & \of & M \\
\end{diagram}

\medskip
The main results of
\cite{FuchsHamkinsReitz:Set-theoreticGeology} attempt to
separate and control these classes.

\begin{theorem}[{(Fuchs, Hamkins, Reitz
\cite{FuchsHamkinsReitz:Set-theoreticGeology})}]
\label{Theorem.MantleCombinations}\
\begin{enumerate}
 \item Every model of set theory $V$ has an extension
     $\Vbar$ with
 $$V=\Mantle^{\Vbar}=\gMantle^{\Vbar}=\gHOD^{\Vbar}=\HOD^{\Vbar}$$
 \item Every model of set theory $V$ has an extension
     $W$ with
 $$V=\Mantle^{W}=\gMantle^{W}=\gHOD^{W}\qquad\text{ but
     }\quad\HOD^{W}=W$$
 \item Every model of set theory $V$ has an extension
     $U$ with
 $$V=\HOD^{U}=\gHOD^{U}\qquad\text{ but }\quad\Mantle^{U}=U$$
 \item Lastly, every $V$ has an extension $Y$ with
 $$Y=\HOD^{Y}=\gHOD^{Y}=\Mantle^{Y}=\gMantle^{Y}$$
\end{enumerate}
\end{theorem}

In particular, as mentioned earlier, every model of \ZFC\
is the mantle and generic mantle of another model of \ZFC.
It follows that we cannot expect to prove any regularity
features about the mantle or the generic mantle, beyond
what can be proved about an arbitrary model of \ZFC. It
also follows that the mantle of $V$ is not necessarily a
ground model of $V$, even when it is a model of \ZFC. One
can therefore iteratively take the mantle of the mantle and
so on, and we have proved that this process can strictly
continue. Indeed, by iteratively computing the mantle of
the mantle and so on, what we call the {\df inner} mantles,
we might eventually arrive at a model of \ZFC\ plus the
ground axiom, where the process would naturally terminate.
If this should occur, then we call this termination model
the {\df outer core} of the original model. Generalizing
the theorem above, Fuchs, Reitz and I have conjectured that
every model of \ZFC\ is the $\alpha^{\rm th}$ inner mantle
of another model of \ZFC, for arbitrary ordinals $\alpha$
(or even $\alpha=\ORD$ or beyond), and we have proved this
for $\alpha<\omega^2$. These inner mantles speak to a
revived version of the ancient paradise position, which
replies to theorems \ref{Theorem.EveryModelIsAMantle} and
\ref{Theorem.MantleCombinations} by saying that if the
mantle turns out itself to be a forcing extension, then one
hasn't really stripped away all the forcing yet, and one
should proceed with the inner mantles, stripping away
additional strata of forcing residue. If the conjecture is
correct, however, then any fixed number of iterations of
this will still arrive essentially at an arbitrary model of
set theory, rather than an ancient set-theoretic paradise.

\bibliographystyle{alpha}
\bibliography{MathBiblio,HamkinsBiblio}

\begin{thebibliography}{FWW08}

\bibitem[Bla84]{Blass1984:InteractionBetweenCategoryTheoryAndSetTheory}
Andreas Blass.
\newblock The interaction between category theory and set theory.
\newblock In {\em Mathematical applications of category theory ({D}enver,
  {C}ol., 1983)}, volume~30 of {\em Contemp. Math.}, pages 5--29. Amer. Math.
  Soc., Providence, RI, 1984.

\bibitem[Coh63]{Cohen1963:IndependenceOfCH}
Paul~J. Cohen.
\newblock The indendependence of the continuum hypothesis.
\newblock {\em Proc. Nat. Acad. Sci. U.S.A.}, 50:1143--1148, 1963.

\bibitem[Coh64]{Cohen1964:IndependenceOfCHII}
Paul~J. Cohen.
\newblock The indendependence of the continuum hypothesis {II}.
\newblock {\em Proc. Nat. Acad. Sci. U.S.A.}, 51:105--110, 1964.

\bibitem[Coh66]{Cohen1966:SetTheoryAndTheContinuumHypothesis}
Paul~J. Cohen.
\newblock {\em Set theory and the continuum hypothesis}.
\newblock W. A. Benjamin, Inc., New York-Amsterdam, 1966.

\bibitem[FHR]{FuchsHamkinsReitz:Set-theoreticGeology}
Gunter Fuchs, Joel~David Hamkins, and Jonas Reitz.
\newblock Set-theoretic geology.
\newblock submitted, preprint available at http://arxiv.org/abs/1107.4776.

\bibitem[Fre86]{Freiling1986:AxiomsOfSymmetry:ThrowingDartsAtTheRealLine}
Chris Freiling.
\newblock Axioms of symmetry: throwing darts at the real number line.
\newblock {\em J. Symbolic Logic}, 51(1):190--200, 1986.

\bibitem[Fri06]{Friedman2006:InternalConsistencyAndIMH}
Sy~D. Friedman.
\newblock Internal consistency and the inner model hypothesis.
\newblock {\em Bulletin of Symbolic Logic}, December 2006.

\bibitem[FWW08]{FriedmanWelchWoodin2008:ConsistencyOfIMH}
Sy-David Friedman, Philip Welch, and W.~Hugh Woodin.
\newblock On the consistency strength of the inner model hypothesis.
\newblock 73:391--400, June 2008.

\bibitem[GH10]{GitmanHamkins2010:NaturalModelOfMultiverseAxioms}
Victoria Gitman and Joel~David Hamkins.
\newblock A natural model of the multiverse axioms.
\newblock {\em Notre Dame Journal of Formal Logic}, 51(4):475--484, 2010.

\bibitem[Hama]{Hamkins:Themultiverse:ANaturalContext}
Joel~David Hamkins.
\newblock The set-theoretical multiverse: a natural context for set theory.
\newblock {\em Annals of the Japan Association for Philosophy of Science,
  special issue Mathematical Logic and Its Applications}.
\newblock to appear.

\bibitem[Hamb]{HamkinsMO6594:AnswerToWhatAreSomeReasonable-soundingStatementsI%
ndepOfZFC?}
Joel~David Hamkins\phantom{x}(mathoverflow.net/users/1946).
\newblock What are some reasonable-sounding statements that are independent of
  {ZFC}?
\newblock MathOverflow.
\newblock \url{http://mathoverflow.net/questions/6594} (version: 2009-11-25).

\bibitem[Ham03a]{Hamkins2003:ExtensionsWithApproximationAndCoverProperties}
Joel~David Hamkins.
\newblock Extensions with the approximation and cover properties have no new
  large cardinals.
\newblock {\em Fundamenta Mathematicae}, 180(3):257--277, 2003.

\bibitem[Ham03b]{Hamkins2003:MaximalityPrinciple}
Joel~David Hamkins.
\newblock A simple maximality principle.
\newblock {\em Journal of Symbolic Logic}, 68(2):527--550, June 2003.

\bibitem[Ham05]{Hamkins2005:TheGroundAxiom}
Joel~David Hamkins.
\newblock The {Ground Axiom}.
\newblock {\em Oberwolfach Report}, 55:3160--3162, 2005.

\bibitem[Ham09]{Hamkins2009:SomeSecondOrderSetTheory}
J.~D. Hamkins.
\newblock Some second order set theory.
\newblock In In~R. Ramanujam and S.~Sarukkai, editors, {\em ICLA 2009, volume
  5378 of LNAI}, pages 36--50. Springer--Verlag, Amsterdam, 2009.

\bibitem[HL08]{HamkinsLoewe2008:TheModalLogicOfForcing}
Joel~David Hamkins and Benedikt L{\"o}we.
\newblock The modal logic of forcing.
\newblock {\em Transactions of the American Mathematical Society},
  360:1793--1817, 2008.

\bibitem[HLR]{HamkinsLinetskyReitz:PointwiseDefinableModelsOfSetTheory}
Joel~David Hamkins, David Linetsky, and Jonas Reitz.
\newblock Pointwise definable models of set theory.
\newblock submitted, preprint available at http://arxiv.org/abs/1105.4597.

\bibitem[HRW08]{HamkinsReitzWoodin2008:TheGroundAxiomAndV=HOD}
Joel~David Hamkins, Jonas Reitz, and W.~Hugh Woodin.
\newblock The {Ground Axiom} is consistent with {V $\neq$ HOD}.
\newblock {\em Proceedings of the AMS}, 136:2943--2949, 2008.

\bibitem[HS]{HamkinsSeabold:BooleanUltrapowers}
Joel~David Hamkins and Daniel Seabold.
\newblock Boolean ultrapowers.
\newblock in preparation.

\bibitem[Kr{\"o}01]{Kromer2001:TarskisAxiomOfInaccessiblesAndGrothendieckUnive%
rses}
Ralf Kr{\"o}mer.
\newblock Tarski's axiom of inaccessibles and {G}rothendieck
  universes---historical and critical remarks on the foundations of category
  theory.
\newblock In {\em Logic, {V}ol. 21 ({P}olish)}, volume 2312 of {\em Acta Univ.
  Wratislav.}, pages 45--57. Wydawn. Uniw. Wroc\l aw., Wroc\l aw, 2001.

\bibitem[Lak67]{Lakatos1967:ProblemsInThePhilosophyOfMathematics}
Imre Lakatos, editor.
\newblock {\em Problems in the philosophy of mathematics}.
\newblock North-Holland Publishing Co., Amsterdam, 1967.
\newblock Proceedings of the International Colloquium in the Philosophy of
  Science, London, 1965, vol. 1.

\bibitem[Lav07]{Laver2007:CertainVeryLargeCardinalsNotCreated}
Richard Laver.
\newblock Certain very large cardinals are not created in small forcing
  extensions.
\newblock {\em Annals of Pure and Applied Logic}, 149(1--3):1--6, 2007.

\bibitem[Luz35]{Luzin1935:SurLesEnssemblesAnalytiquesNuls}
N.~N. Luzin.
\newblock Sur les ensembles analytiques nuls.
\newblock {\em Fund. Math.}, 25:109--131, 1935.

\bibitem[Mad98]{Maddy1998:V=LAndMaximize}
Penelope Maddy.
\newblock {$V=L$} and {Maximize}.
\newblock {\em Logic Colloquium '95 Proceedings of ASL, Haifa, Isreal 1995},
  pages 134--52, 1998.
\newblock Makowski J. A. and Raave E. V., editors.

\bibitem[Mar01]{Martin2001:MultipleUniversesOfSetsAndIndeterminateTruthValues}
Donald~A. Martin.
\newblock Multiple universes of sets and indeterminate truth values.
\newblock {\em Topoi}, 20(1):5--16, 2001.

\bibitem[Pea89]{Peano1889:PrinciplesOfArithmeticPresentedByNewMethod}
Giuseppe Peano.
\newblock The principles of arithmetic, presented by a new method.
\newblock pages pp. 83--97 in van Heijenoort, 1967, 1889.

\bibitem[Rei06]{Reitz2006:Dissertation}
Jonas Reitz.
\newblock {\em The {Ground Axiom}}.
\newblock PhD thesis, The Graduate Center of the City University of New York,
  September 2006.

\bibitem[Rei07]{Reitz2007:TheGroundAxiom}
Jonas Reitz.
\newblock The {Ground Axiom}.
\newblock {\em Journal of Symbolic Logic}, 72(4):1299--1317, 2007.

\bibitem[Sol]{Solovay2010:PersonalCommunication}
Robert~M. Solovay.
\newblock personal communication, October 31, 2010.

\bibitem[Sol76]{Solovay1976:ProvabilityModalLogic}
Robert~M. Solovay.
\newblock Provability interpretations of modal logic.
\newblock {\em Israel Journal of Mathematics}, 25:287--304, 1976.

\bibitem[Ste]{Steel2004:SlidesGenericAbsolutenessAndTheContinuumProblem}
John~R. Steel.
\newblock Generic absoluteness and the continuum problem.
\newblock Slides for a talk at the Laguna workshop on philosophy and the
  continuum problem (P. Maddy and D. Malament organizers) March, 2004.
\newblock http://math.berkeley.edu/~steel/talks/laguna.ps.

\bibitem[SV03]{StaviVaananen:ReflectionPrinciples}
J.~Stavi and J.~V\"a\"an\"anen.
\newblock Reflection principles for the continuum.
\newblock In Yi~Zhang, editor, {\em Logic and Algebra}, volume 302 of {\em AMS
  Contemporary Mathematics Series}. 2003.

\bibitem[vH67]{vanHeijenoort1967:FromFregeToGodelSourcebookInMathLogic}
Jean van Heijenoort.
\newblock {\em From {F}rege to {G}\"odel. {A} source book in mathematical
  logic, 1879--1931}.
\newblock Harvard University Press, Cambridge, Mass., 1967.

\bibitem[Woo04]{Woodin2004:RecentDevelopmentsOnCH}
W.~Hugh Woodin.
\newblock Recent development's on {Cantor's Continuum Hypothesis}.
\newblock {\em Proceedings of the Continuum in Philosophy and Mathematics},
  2004.
\newblock Carlsberg Academy, Copenhagen, November 2004.

\bibitem[Zer30]{Zermelo1930:UberGrenzzahlenUndMengenBereiche}
E.~Zermelo.
\newblock {\"U}ber {G}renzzahlen und {M}engenbereiche. {N}eue {U}ntersuchungen
  \"uber die {G}rundlagen der {M}engenlehre. ({German}).
\newblock {\em Fundumenta Math.}, 16:29--47, 1930.

\end{thebibliography}

\end{document}